\documentclass{tran-l}
\usepackage[mathscr]{eucal}
\usepackage{amssymb}
\usepackage{enumerate}
\usepackage{graphpap}
\usepackage{mathrsfs}

\newtheorem{theorem}{Theorem}[section]
\newtheorem{corollary}[theorem]{Corollary}
\newtheorem{lemma}[theorem]{Lemma}
\newtheorem{proposition}[theorem]{Proposition}

\theoremstyle{definition}

\newtheorem{example}[theorem]{Example}
\newtheorem{question}[theorem]{Question}

\theoremstyle{remark}
\newtheorem{remark}[theorem]{Remark}

\numberwithin{equation}{section}

\newcommand{\up}{\textup}
\newcommand{\con}{\operatorname{\mathsf{con}}}
\newcommand{\ini}{\operatorname{\mathsf{ini}}}

\newcommand{\occ}{\operatorname{\mathsf{occ}}}
\newcommand{\ba}{\mathbf{a}}
\newcommand{\bb}{\mathbf{b}}
\newcommand{\bc}{\mathbf{c}}
\newcommand{\bd}{\mathbf{d}}
\newcommand{\be}{\mathbf{e}}
\newcommand{\bH}{\mathbf{H}}
\newcommand{\bp}{\mathbf{p}}
\newcommand{\bq}{\mathbf{q}}
\newcommand{\bu}{\mathbf{u}}
\newcommand{\bV}{\mathbf{V}}
\newcommand{\bv}{\mathbf{v}}
\newcommand{\bw}{\mathbf{w}}
\newcommand{\bz}{\mathbf{z}}

\newcommand{\bbB}{\mathbb{B}}
\newcommand{\bbE}{\mathbb{E}}
\newcommand{\bbI}{\mathbb{I}}
\newcommand{\bbJ}{\mathbb{J}}
\newcommand{\bbL}{\mathbb{L}}
\newcommand{\bbM}{\mathbb{M}}
\newcommand{\bbQ}{\mathbb{Q}}
\newcommand{\bbR}{\mathbb{R}}
\newcommand{\bbV}{\mathbb{V}}

\newcommand{\cK}{\mathcal{K}}
\newcommand{\cM}{\mathcal{M}}
\newcommand{\cU}{\mathcal{U}}
\newcommand{\cV}{\mathcal{V}}

\newcommand{\fL}{\mathfrak{L}}

\newcommand{\one}{type~$\aleph_0$}
\newcommand{\two}{type~$2^{\aleph_0}$}

\newcommand{\RQ}{\mathsf{M}}

\newcommand{\sW}{\mathscr{W}}
\newcommand{\sA}{\mathscr{A}}
\newcommand{\sX}{\mathscr{X}}

\newcommand{\Vm}{\mathbb{V}_\mathrm{mon}}
\newcommand{\Vs}{\mathbf{V}_\mathrm{sem}}

\newcommand{\hfb}{hereditarily finitely based}
\newcommand{\infb}{in\-her\-ently \nfb}
\newcommand{\nfb}{non\-finitely based}

\allowdisplaybreaks

\begin{document}

\title{Monoid varieties with extreme properties}

\author[M. Jackson]{Marcel Jackson}
\address{Department of Mathematics and Statistics, La Trobe University, Victoria~3086, Australia}
\email{m.g.jackson@latrobe.edu.au}
\thanks{The first author was supported by ARC Discovery Project DP1094578 and Future Fellowship FT120100666}

\author[E. W. H. Lee]{Edmond W. H. Lee}
\address{Department of Mathematics, Nova Southeastern University, Florida~33314, USA}
\email{edmond.lee@nova.edu}


\dedicatory{Dedicated to the 81st birthday of John L. Rhodes}

\begin{abstract}
Finite monoids that generate monoid varieties with uncountably many subvarieties seem rare, and surprisingly, no finite monoid is known to generate a monoid variety with countably infinitely many subvarieties.
In the present article, it is shown that there are, nevertheless, many finite monoids with simple descriptions that generate monoid varieties with continuum many subvarieties; these include inherently non\-finitely based finite monoids and all monoids for which $xyxy$ is an isoterm.
It follows that the join of two Cross monoid varieties can have a continuum cardinality subvariety lattice that violates the ascending chain condition.

Regarding monoid varieties with countably infinitely many subvarieties, the first example of a finite monoid that generates such a variety is exhibited.
A complete description of the subvariety lattice of this variety is given.
This lattice has width three and contains only finitely based varieties, all except two of which are Cross.
\end{abstract}

\maketitle

\section{Introduction}

A \textit{monoid} is a semigroup with an identity element.
Any semigroup~$S$ without identity element can be converted into the monoid ${S^1 = S \cup \{1\}}$, where multiplication involving~$1$ is defined by ${1 \cdot x = x \cdot 1 = x}$ for all ${x \in S^1}$.
In the present article, the notion of a \textit{monoid} is attached ambiguously to either a semigroup with identity element in the plain signature ${ \{\,\cdot\,\} }$ or to a semigroup with identity element~$1$ in the enriched signature $\{\,\cdot\,,1\}$.
To distinguish this difference for any class~$\cM$ of monoids at the varietal level, let $\Vs\cM$ denote semigroup variety generated by~$\cM$ and let $\Vm\cM$ denote the monoid variety generated by~$\cM$.
More generally, bold upper case letters ${ \mathbf{A}, \mathbf{B}, \mathbf{C}, \ldots }$ represent semigroup varieties, while blackboard bold upper case letters ${ \mathbb{A}, \mathbb{B}, \mathbb{C}, \ldots }$ represent monoid varieties.

In the semigroup signature ${ \{\,\cdot\,\} }$, monoids have played a conspicuous role in the study of semigroup varieties with extreme properties.
For instance, a finite semigroup is {\infb} if and only if it contains an {\infb} submonoid~\cite{mS88a,mS88b}; recall that an algebra~$A$ is \textit{\infb} if every locally finite variety containing~$A$ is {\nfb}.
The first two published examples of {\nfb} finite semigroups, due to Perkins~\cite{Per69}, are monoids.
The first example---the monoid~$B_2^1$ obtained from the Brandt semigroup
\[
B_2 = \big\langle\, a,b \,\big|\, a^2=b^2=0, \, aba=a, \, bab=b \,\big\rangle
\]
of order five---is now known to be one of only four {\nfb} semigroups or order six~\cite{LLZ12}; the other example is a monoid that is \textit{nilpotent} in the sense that every sufficiently long product not involving the identity element is equal to the zero element.
Nilpotent monoids emerged as a significant source of {\nfb} finite semigroups~\cite{Jac01,Jac05a,JS00,oS00,oS15} and of semigroups that generate varieties with continuum many subvarieties~\cite{Jac00}.

A variety of algebras that contains finitely many subvarieties is said to be \textit{small}.
A \textit{Cross variety} is a variety that is finitely based, finitely generated, and small.
The present article is concerned with monoid varieties that are not small.
A variety with countably infinitely many subvarieties is said to be of~\textit{\one}, while a variety with continuum many subvarieties is said to be of~\textit{\two}.

\subsection{Finite monoids that generate monoid varieties of~\two} \label{subsec: continuum}

The subvarieties of a variety~$\cV$ of algebras constitute a complete lattice $\fL(\cV)$ with respect to class inclusion.
For any monoid~$M$, the lattice $\fL(\Vm\{M\})$ order-embeds into the lattice $\fL(\Vs\{M\})$ in the obvious manner, so in one sense, the latter lattice is as complicated as the former.
However, there exist examples of finite monoids~$M$ where the monoid variety $\Vm\{M\}$ is small even though the semigroup variety $\Vs\{M\}$ is of~\two; the monoid~$N_6^1$ obtained from the nilpotent semigroup
\[
N_6 = \big\langle\, a,b \,\big|\, a^2=b^2=aba=0 \,\big\rangle
\]
of order six is one such example~\cite{Jac00,Jac05b}.
At the moment, the only examples of finite monoids that generate monoid varieties of~{\two} are those of Jackson and McKenzie~\cite{JM06} arising from certain graph-encoding technique.
But these monoids, the smallest of which is of order~20, have a somewhat nontrivial description.
In the present article, results of Jackson~\cite{Jac00} are extended to the monoid setting to exhibit new easily described examples of monoids that generate monoid varieties of~\two.
In particular, these examples include many nilpotent monoids and all {\infb} finite monoids.
It is also shown that the join of two Cross monoid varieties can be of~\two.

\subsection{Finite monoids that generate varieties of~\one}

Many finite semigroups are known to generate varieties of~\one.
For instance, an abundance of examples can be found from varieties generated by completely 0-simple semigroups~\cite{Lee10,LV11}; the Brandt semigroup~$B_2$ is an example~\cite{Lee06} in particular.

Another easy method of locating varieties of~{\one} is to consider \textit{\hfb} varieties, that is, varieties all subvarieties of which are finitely based.
Since every {\hfb} variety contains at most countably many subvarieties, non\-small {\hfb} varieties are of~\one.
For example, commutative semigroups are long known to be \hfb~\cite{Per69}; since the commutative monoid~$N_2^1$ obtained from the nilpotent semigroup
\[
N_2 = \big\langle\, a \,\big|\, a^2=0 \,\big\rangle
\]
generates a non\-small semigroup variety~\cite[Figure~5(b)]{Eva71}, the direct product of~$N_2^1$ with any finite commutative semigroup generates a semigroup variety of~\one.

The situation changes drastically when monoid varieties are considered.
It so happens that {\hfb} monoid varieties generated by presently known finite monoids---commutative monoids~\cite{Hea68}, idempotent monoids~\cite{Bir70,Fen71,Ger70}, nilpotent monoids~\cite{Lee11}, and 2-testable monoids~\cite{Lee12} to name a few---are all Cross and so are not of~\one.
In fact, no example of a finite monoid is known to generate a monoid variety of~\one.
The present article exhibits the first such example.
Specifically, a certain monoid~$E^1$ of order six is shown to generate a non\-small, {\hfb} monoid variety.

\subsection{Organization}

There are six sections in the present article.
Notation and background material are given in~\S\ref{sec: prelim}.
In~\S\ref{sec: xyxy}, a certain nilpotent monoid of order nine is shown to generate a monoid variety of~\two.
Based on this result, a number of other monoid varieties of~{\two} are also exhibited.
In~\S\ref{sec: INFB}, all inherently {\nfb} finite monoids are shown to generate monoid varieties of~\two.
The two well-known {\infb} monoids of order six are thus small examples that generate monoid varieties of~\two; whether or not these two examples are minimal is discussed.
In~\S\ref{sec: E}, the subvariety lattice of a monoid variety generated by a certain monoid~$E^1$ of order six is completely described; in particular, this variety is non\-small and {\hfb} and so is of~\one.
In~\S\ref{sec: other}, some other extreme properties satisfied by monoid varieties are discussed and related open questions are posed.

\section{Preliminaries} \label{sec: prelim}

\subsection{Words and identities}

Let~$\sA$ be a countably infinite alphabet.
For any subset~$\sX$ of~$\sA$, let~$\sX^+$ and~$\sX^*$ denote the free semigroup and free monoid over~$\sX$, respectively.
Elements of~$\sA$ are called \textit{variables}, and elements of~$\sA^*$ are called \textit{words}.
The empty word, more conveniently written as~$1$, is the identity element of the monoid~$\sA^*$.
For any word~$\bw$,
\begin{itemize}
\item the \textit{content} of~$\bw$, denoted by $\con(\bw)$, is the set of variables occurring in~$\bw$;
\item the number of occurrences of a variable~$x$ in~$\bw$ is denoted by $\occ(x,\bw)$;
\item a variable~$x$ in~$\bw$ is \textit{simple} if $\occ(x,\bw)=1$;
\item the \textit{initial part} of~$\bw$, denoted by $\ini(\bw)$, is the word obtained by retaining the first occurrence of each variable in~$\bw$.
\end{itemize}
For any ${ x_1,x_2,\ldots,x_r \in \sA }$, let ${ \bw[x_1,x_2,\ldots,x_r] }$ denote the word obtained by applying the substitution that fixes ${ x_1,x_2,\ldots,x_r }$ and assigns the value~$1$ to all other variables.
For example, ${ \bw[x]=x^{\occ(x,\bw)} }$ for any ${ x\in \con(\bw) }$.

An \textit{identity} is an expression ${ \bu \approx \bv }$, where ${ \bu,\bv \in \sA^* }$.
A semigroup~$S$ \textit{satisfies} an identity ${ \bu \approx \bv }$, written ${ S \models \bu \approx \bv }$, if for any substitution ${\theta: \sA \to S}$, the equality ${ \bu\theta=\bv\theta }$ holds in~$S$.
For any class~$\cK$ of semigroups, write ${ \cK \models \bu \approx \bv }$ to indicate that ${ S \models \bu \approx \bv }$ for all ${ S \in \cK }$.
A monoid that satisfies an identity ${ \bu \approx \bv }$ clearly also satisfies the identity ${ \bu[x_1,x_2,\ldots,x_r] \approx \bv[x_1,x_2,\ldots,x_r] }$ for any ${ x_1,x_2,\ldots,x_r \in \sA }$.

An identity ${ \bw \approx \bw' }$ is \textit{directly deducible} from an identity ${ \bu \approx \bv }$ if some substitution ${ \theta: \sA \to \sA^* }$ and words ${ \ba, \bb \in \sA^* }$ exist such that ${ \{\bw, \bw'\} = \{\ba(\bu\theta)\bb,\ba(\bv\theta)\bb\} }$.
An identity ${ \bw \approx \bw' }$ is \textit{deducible} from a set~$\Sigma$ of identities, written ${ \Sigma \vdash \bw \approx \bw' }$, if there exists a sequence
\[
\bw =\bu_0, \bu_1, \ldots, \bu_r = \bw'
\]
of distinct words such that each identity ${ \bu_i \approx \bu_{i+1} }$ is directly deducible from some identity in~$\Sigma$.
By Birkhoff's completeness theorem of equational logic, a deduction ${ \Sigma \vdash \bw \approx \bw' }$ holds if and only if any monoid that satisfies the identities in~$\Sigma$ also satisfies the identity ${ \bw \approx \bw' }$; see Burris and Sankappanavar~\cite[Theorem~14.19]{BS81}.

\subsection{Varieties}

A variety is a class of algebras that is closed under the formation of homomorphic images, subalgebras, and arbitrary direct products.
Equivalently, a variety~$\cV$ is a class of algebras that satisfy some set~$\Sigma$ of identities; in this case, $\Sigma$ is said to be a \textit{basis} for~$\cV$, or~$\cV$ is \textit{defined} by~$\Sigma$.
For any variety~$\cV$ and set~$\Sigma$ of identities, let $\cV\Sigma$ denote the subvariety of~$\cV$ \textit{defined} by~$\Sigma$.

A variety is \textit{finitely based} if it has some finite basis, \textit{finitely generated} if it is generated by some finite algebra, and \textit{small} if it contains finitely many subvarieties.
A variety is \textit{Cross} if it is finitely based, finitely generated, and small.

The following well- and long-known result can be found, after easy modification, in Oates-MacDonald and Vaughan-Lee \cite[p.~370]{OV78}.

\begin{lemma} \label{L: FG}
A locally finite variety~$\cV$ of algebras is finitely generated if and only if there exists no strictly increasing infinite chain ${ \cV_1 \subset \cV_2 \subset \cV_3 \subset \cdots }$ of varieties such that ${ \cV = \bigvee_{i \geq 1} \cV_i }$\up.
Consequently\up, any locally finite\up, small variety is finitely generated\up.
\end{lemma}

For any subvariety~$\cU$ of a variety~$\cV$, let ${ [\cU,\cV] }$ denote the interval of all subvarieties of~$\cV$ containing~$\cU$.
Let $\fL(\cV)$ denote the lattice of all subvarieties of~$\cV$, that is, ${ \fL(\cV) = [{\bf 0},\cV] }$, where~$\bf 0$ is the trivial variety.

For any class~$\cM$ of monoids, let $\Vs\cM$ denote semigroup variety generated by~$\cM$ and let $\Vm\cM$ denote the monoid variety generated by~$\cM$.

\begin{lemma}[Jackson~{\cite[Lemma~1.1]{Jac05c}}] \label{L: S S1}
Let~$\cM$ be any class of monoids and let~$S$ be any semigroup\up.
Then ${ S \in \Vs\cM }$ implies that ${ S^1 \in \Vm\cM }$\up.
\end{lemma}

\subsection{Rees quotients of free monoids}

For any words ${ \bu,\bv \in \sA^* }$, the expression ${ \bu \preceq \bv }$ indicates that~$\bu$ is a \textit{factor} of~$\bv$, that is, ${ \bv \in \sA^* \bu \sA^* }$.
For any set ${ \sW \subseteq \sA^* }$ of words, define the \textit{factorial closure} of~$\sW$ to be the set~$\sW^\preceq$ of all words that are factors of some word in~$\sW$:
\[
\sW^\preceq=\{\bu \in \sA^* \,|\, \text{$\bu \preceq \bw$ for some $\bw \in \sW$}\}.
\]
Then ${ I_\sW = \sA^* \backslash \sW^\preceq }$ is an ideal of~$\sA^*$.
Let~$\RQ\sW$ denote the Rees quotient $\sA^*/I_\sW$.
Equivalently, $\RQ\sW$ can be treated as the monoid that consists of every factor of every word in~$\sW$, together with a zero element~$0$, with binary operation~$\cdot$ given by
\[
\bu \cdot \bv =
\begin{cases}
\,\bu\bv & \text{if $\bu\bv$ is a factor of some word in~$\sW$,} \\
\,0 & \text{otherwise}.
\end{cases}
\]

\begin{example}
The monoids $\RQ\{x\}$ and~$\RQ\{xyx\}$ are isomorphic to the monoids~$N_2^1$ and~$N_6^1$ in~\S\ref{subsec: continuum}, respectively, while ${ \RQ\varnothing = \{ 0,1\} }$ is the semilattice of order two.
\end{example}

\begin{remark}
The Rees quotient ${ \sA^*/I_\sW }$ is more commonly denoted in the literature by $\mathrm{S}(\sW)$; see, for example, Jackson~\cite{Jac00,Jac01,Jac05b}, Jackson and~Sapir~\cite{JS00}, Lee~\cite{Lee11,Lee12,Lee14a}, and O.~Sapir~\cite{oS00}.
But in retrospect, the symbol $\mathrm{S}(\sW)$ might have been better used to denote the Rees quotient $\sA^+/I_\sW$ without identity element.
\end{remark}

A word~$\bw$ is an \textit{isoterm} for a class~$\cM$ of monoids if whenever ${ \cM \models \bw \approx \bw' }$ for some ${ \bw' \in \sA^+ }$, then ${ \bw \approx \bw' }$ is satisfied by all monoids; in other words, once all instances of~$1$ have been removed by applications of ${ x1 \mapsto x }$ and ${ 1x \mapsto x }$, the two words~$\bw$ and~$\bw'$ are identical.

\begin{lemma}[Jackson~{\cite[Lemma~3.3]{Jac05b}}] \label{L: isoterm}
For any class~$\cM$ of monoids and any set~$\sW$ of words\up,
${ \RQ\sW\in \Vm\cM }$ if and only if each word in~$\sW$ is an isoterm for~$\cM$\up.
\end{lemma}

\subsection{Varieties generated by some small monoids}

The monoids~$B_0^1$, $I^1$, $J^1$, $L_2^1$, and~$R_2^1$ obtained from the following semigroups are required in later sections:
\begin{align*}
B_0 & = \big\langle a,b,c \,\big|\, a^2 = a, \, b^2 = b, \, ab=ba=0, \, ac = cb = c \big\rangle, \\
I   & = \big\langle a,b \,\big|\, ab = a, \, ba = 0, \, b^2 = b \big\rangle, \\
J   & = \big\langle a,b \,\big|\, ba = a, \, ab = 0, \, b^2 = b \big\rangle, \\
L_2 & = \big\langle a,b \,\big|\, a^2 = ab = a, \, b^2 = ba = b \big\rangle, \\
R_2 & = \big\langle a,b \,\big|\, a^2 = ba = a, \, b^2 = ab = b \big\rangle.
\end{align*}
Note that~$B_0$, $I$, and~$J$ are subsemigroups of the Brandt semigroup~$B_2$, while~$L_2$ and~$R_2$ are the left- and right-zero semigroups of order two, respectively.
Write
\begin{gather*}
\bbB_0^1 = \Vm\{B_0^1\}, \quad \bbI^1 = \Vm\{I^1\}, \quad \bbJ^1 = \Vm\{J^1\}, \\ \bbL_2^1 = \Vm\{L_2^1\}, \quad \text{and} \quad \bbR_2^1 = \Vm\{R_2^1\}.
\end{gather*}
For any set ${ \sW \subseteq \sA^* }$ of words, write
\[
\bbM\sW = \Vm \{\RQ\sW\}.
\]

\begin{proposition}[Lee~{\cite[Proposition~4.1]{Lee12}}] \label{P: LNR}
The monoid variety ${ \bbL_2^1 \vee \bbM\{x\} \vee \bbR_2^1 }$ is Cross\up, and the lattice ${ \fL(\bbL_2^1 \vee \bbM\{x\} \vee \bbR_2^1) }$ is given in Figure~\ref{F: LNR}\up.
\end{proposition}

\begin{figure}[htbp]
\begin{picture}(290,160)(0,0) \setlength{\unitlength}{0.50mm}
\put(100,100){\circle*{2}}
\put(102,103){\makebox(0,0)[b]{$\bbL_2^1 \vee \bbM\{x\} \vee \bbR_2^1$}}
\put(100,80){\circle*{2}}
\put(103,79){\makebox(0,0)[bl]{$\bbL_2^1 \vee \bbR_2^1$}}
\put(100,80){\line(0,1){20}}
\put(60,80){\circle*{2}}
\put(58,82){\makebox(0,0)[br]{$\bbL_2^1 \vee \bbB_0^1$}}
\put(40,70){\circle*{2}}
\put(36,70){\makebox(0,0)[r]{$\bbL_2^1 \vee \bbM\{x\}$}}
\put(40,50){\circle*{2}}
\put(36,50){\makebox(0,0)[r]{$\bbL_2^1$}}
\put(40,50){\line(0,1){20}}
\put(140,80){\circle*{2}}
\put(142,81){\makebox(0,0)[bl]{$\bbB_0^1 \vee \bbR_2^1$}}
\put(160,70){\circle*{2}}
\put(165,70){\makebox(0,0)[l]{$\bbM\{x\} \vee \bbR_2^1$}}
\put(160,50){\circle*{2}}
\put(165,50){\makebox(0,0)[l]{$\bbR_2^1$}}
\put(160,50){\line(0,1){20}}
\put(100,60){\circle*{2}}
\put(100,64){\makebox(0,0)[b]{$\bbB_0^1$}}
\put(80,50){\circle*{2}}
\put(77,51){\makebox(0,0)[tr]{$\bbI^1$}}
\put(120,50){\circle*{2}}
\put(124,51){\makebox(0,0)[tl]{$\bbJ^1$}}
\put(100,40){\circle*{2}}
\put(104,41){\makebox(0,0)[tl]{$\bbM\{xy\}$}}
\put(100,30){\circle*{2}}
\put(98,29){\makebox(0,0)[br]{$\bbM\{x\}$}}
\put(100,20){\circle*{2}}
\put(104,21){\makebox(0,0)[tl]{$\bbM\varnothing$}}
\put(100,10){\circle*{2}}
\put(100,05){\makebox(0,0)[t]{$\mathbf{0}$}}
\put(100,40){\line(0,-1){30}}
\put(100,100){\line(-2,-1){60}}\put(100,100){\line(2,-1){60}}
\put(80,50){\line(2,1){60}}\put(120,50){\line(-2,1){60}}
\put(100,40){\line(-2,1){60}}\put(100,40){\line(2,1){60}}
\put(100,80){\line(-2,-1){60}}\put(100,80){\line(2,-1){60}}
\put(100,20){\line(2,1){60}}\put(100,20){\line(-2,1){60}}
%
\end{picture}\caption{The lattice $\fL(\bbL_2^1 \vee \bbM\{x\} \vee \bbR_2^1)$}
\label{F: LNR}
\end{figure}

\subsection{Words with overlapping variables} \label{subsec: overlapping}

The monoid $\RQ\{xyx\}$ plays a crucial role in the investigation of semigroup varieties of~\two.

\begin{proposition}[Jackson~{\cite[\S3]{Jac00}}] \label{P: xyx}
The variety $\Vs\big\{\RQ\{xyx\}\big\}$ is of~\two\up.
\end{proposition}

In the proof of Proposition~\ref{P: xyx}, Jackson established the irredundancy of the following identity system within the equational theory of $\Vs\big\{\RQ\{xyx\}\big\}$:
\begin{equation}
\{ x\bu_n \approx \bu_n, \, \bu_n x \approx \bu_n \,|\, n = 3,4,5,\ldots \} \label{system}
\end{equation}
where
\[
\bu_n = hx_1x_2x_3x_4h \cdot y_1z_1y_1 \cdot y_2z_2y_2 \cdots y_nz_ny_n \cdot tx_5x_6x_7x_8t.
\]
More specifically, Jackson~\cite[proof of Theorem~3.2]{Jac00} demonstrated that for any identity basis~$\Sigma$ for $\RQ\{xyx\}$ and any subset~$\mathtt{N}$ of ${ \{3,4,5,\ldots\} }$, the deduction
\[
\Sigma \cup \{ x\bu_n \approx \bu_n, \, \bu_n x \approx \bu_n \,|\, n \in \mathtt{N} \} \vdash \{ x\bu_k \approx \bu_k, \, \bu_k x \approx \bu_k\}
\]
holds if and only if ${ k \in \mathtt{N} }$.
It follows that each subset of ${ \{3,4,5,\ldots\} }$ corresponds uniquely to a subvariety of $\Vs\big\{\RQ\{xyx\}\big\}$, thus establishing Proposition~\ref{P: xyx}.

One main goal of the present article is to exhibit finite monoids that generate monoid varieties of~\two.
The monoid $\RQ\{xyx\}$ does not serve this purpose because it generates a monoid variety with only five subvarieties~\cite[Lemma~4.4]{Jac05b}.
It turns out that words of the following form with overlapping variables will be useful in all of the examples of the present article regarding monoid varieties of~\two:
\begin{equation}
x_0 \ \underline{\,?\,} \ x_1 \ \underline{\,?\,} \ x_0 \cdot x_2x_1 \cdot x_3x_2 \cdot x_4x_3\cdots x_{n-1} x_{n-2} \cdot x_n \ \underline{\,?\,} \ x_{n-1} \ \underline{\,?\,} \ x_n \label{pattern}
\end{equation}
The segments indicated by question marks are either empty or filled by variables different from ${ x_0,x_1,\ldots,x_n }$.
Note that crucially, assigning~$1$ to any of the variables ${ x_0,x_1,\ldots,x_n }$ destroys the overlapping pattern; this eliminates applicability of such reduced words to similar patterns in different numbers of variables.
In particular, ``short" identities formed by words of the form~\eqref{pattern} are not deducible from ``long" identities similarly formed.

\begin{remark}
Identities formed by words of the form~\eqref{pattern} were first employed by Lee~\cite{Lee14a,Lee14b} and Lee and Zhang~\cite{LZ14} to exhibit several examples of non\-finitely generated monoid varieties.
\end{remark}

\section{New examples of monoid varieties of~\two} \label{sec: xyxy}

The main example of this section is the monoid variety~$\bbM\{xyxy\}$.
In~\S\ref{subsec: xyxy}, this variety is shown to be of~\two.
Therefore any monoid variety that contains the monoid~$\RQ\{xyxy\}$ is also of~\two.
In particular, it is shown in~\S\ref{subsec: join} that the join of two Cross monoid varieties can be of~\two.
Two open questions regarding semigroup varieties are also shown to have answers within the context of monoid varieties.
In~\S\ref{subsec: order n}, the monoid variety generated by all monoids of order~$n$ is investigated.
It is shown that this variety is of~{\two} if and only if $n \geq 4$.

\subsection{The variety $\bbM\{xyxy\}$} \label{subsec: xyxy}

In the present subsection, define the words
\begin{align*}
\bw_n & = x_0 \cdot yz \cdot x_1x_0 \cdot x_2x_1 \cdot x_3x_2 \cdots x_nx_{n-1} \cdot yz \cdot x_n \\
\text{and} \quad \bw_n' & = x_0 \cdot zy \cdot x_1x_0 \cdot x_2x_1 \cdot x_3x_2 \cdots x_nx_{n-1} \cdot zy \cdot x_n.
\end{align*}

\begin{lemma} \label{L: isoterms for S(xyxy)}
For each ${ n \geq 3 }$\up, the words~$\bw_n$ and~$\bw_n'$ are isoterms for $\RQ\{xyxy\}$\up.
\end{lemma}

\begin{proof}
By symmetry, it suffices to show that~$\bw_n$ is an isoterm for $\RQ\{xyxy\}$.
Let ${ \bw_n \approx \bu }$ be any identity satisfied by $\RQ\{xyxy\}$.
Then the identities
\begin{gather*}
\bw_n[y,z] \approx \bu[y,z], \quad \bw_n[y,x_0] \approx \bu[y,x_0], \quad \bw_n[y,x_n] \approx \bu[y,x_n], \\
\bw_n[z,x_0] \approx \bu[z,x_0], \quad \bw_n[z,x_n] \approx \bu[z,x_n], \quad \text{and} \quad \bw_n[x_i,x_{i+1}] \approx \bu[x_i,x_{i+1}]
\end{gather*}
are satisfied by $\RQ\{xyxy\}$.
Since
\begin{gather*}
\bw_n[y,z] =yzyz, \quad \bw_n[y,x_0] =x_0y x_0y, \quad \bw_n[y,x_n] =y x_n y x_n, \\
\bw_n[z,x_0] = x_0z x_0z, \quad \bw_n[z,x_n] = z x_nz x_n, \quad \text{and} \quad \bw_n[x_i,x_{i+1}] = x_ix_{i+1}x_ix_{i+1}
\end{gather*}
are isoterms for $\RQ\{xyxy\}$, it follows that
\begin{enumerate}
\item[(a)] ${ \bu[x_i,x_{i+1}] = x_ix_{i+1}x_ix_{i+1} }$ for all ${ i \in \{0,1,\ldots, n-1\} }$,
\item[(b)] ${ \bu[y,x_0]=x_0y x_0y }$, \ ${ \bu[y,x_n] =y x_ny x_n }$,
\item[(c)] ${ \bu[z,x_0] = x_0z x_0z }$, \ ${ \bu[z,x_n] = z x_nz x_n }$,
\item[(d)] ${ \bu[y,z]=yzyz }$.
\end{enumerate}
Then~(a) implies that for any ${ i \in \{0,1,\ldots, n-2\} }$, the word ${ \bu[x_i,x_{i+1},x_{i+2}] }$ is
\[
\bp = x_ix_{i+1}x_ix_{i+2}x_{i+1}x_{i+2} \quad \text{or} \quad \bq = x_ix_{i+1}x_{i+2}x_ix_{i+1}x_{i+2}.
\]
Suppose that ${ \bu[x_i,x_{i+1},x_{i+2}] = \bq }$.
Then since ${ \bw_n[x_i,x_{i+1},x_{i+2}] \approx \bu[x_i,x_{i+1},x_{i+2}] }$ is an identity of~$\RQ\{xyxy\}$ and ${ \bw_n[x_i,x_{i+1},x_{i+2}] = \bp }$, it follows that ${ \bp \approx \bq }$ is an identity of~$\RQ\{xyxy\}$.
Therefore ${ \bp[x_i,x_{i+2}] \approx \bq[x_i,x_{i+2}] }$ is also an identity of~$\RQ\{xyxy\}$; but this is impossible because ${ \bq [x_i,x_{i+2}] = x_ix_{i+2}x_ix_{i+2} }$ is an isoterm for~$\RQ\{xyxy\}$, while ${ \bp[x_i,x_{i+2}] = x_i^2x_{i+2}^2 }$.
Hence ${ \bu[x_i,x_{i+1},x_{i+2}] = \bp }$, that is,
\begin{enumerate}
\item[(e)] ${ \bu[x_i,x_{i+1},x_{i+2}] = x_ix_{i+1}x_ix_{i+2}x_{i+1}x_{i+2} }$ for all ${ i \in \{0,1,\ldots, n-2\} }$.
\end{enumerate}
It is then routinely shown by~(a) and~(e) that
\begin{enumerate}
\item[(f)] ${ \bu[x_0,x_1,\ldots,x_n] = x_0 \cdot x_1x_0 \cdot x_2x_1 \cdot x_3x_2 \cdots x_{n-1}x_{n-2} \cdot x_nx_{n-1} \cdot x_n }$.
\end{enumerate}

Now~(b) implies that ${ \occ(y,\bu)=2 }$, where the first~$y$ of~$\bu$ is sandwiched between the two occurrences of~$x_0$, while the second~$y$ of~$\bu$ is sandwiched between the two occurrences of~$x_n$.
Hence by~(f),
\begin{align*}
\bu[y,x_0,x_1,\ldots,x_n] = x_0 \cdot y^{r_1} x_1 y^{r_2} x_0 \cdot x_2x_1 \cdot x_3x_2 \cdots x_{n-1}x_{n-2} \cdot x_n y^{s_1} x_{n-1} y^{s_2} \cdot x_n
\end{align*}
for some ${ (r_1,r_2), (s_1,s_2) \in \{(1,0),(0,1)\} }$.
Note that ${ \bw_n[y,x_1] \approx \bu[y,x_1] }$ is an iden\-tity of~$\RQ\{xyxy\}$, where ${ \bw_n[y,x_1] = y x_1^2y }$ and ${ \bu[y,x_1] = y^{r_1} x_1 y^{r_2} x_1 y }$.
Since the word $x_1y x_1y$ is an isoterm for~$\RQ\{xyxy\}$, it follows that ${ \bu[y,x_1] \neq x_1y x_1y }$ so that ${ (r_1,r_2) \neq (0,1) }$.
By a symmetrical argument, ${ \bw_n[y,x_{n-1}] \approx \bu[y,x_{n-1}] }$ is an identity of $\RQ\{xyxy\}$, where ${ \bw_n[y,x_{n-1}] = y x_{n-1}^2y }$ and ${ \bu[y,x_{n-1}] = y x_{n-1}y^{s_1}x_{n-1}y^{s_2} }$, whence ${ (s_1,s_2) \neq (1,0) }$.
Therefore ${ (r_1,r_2) = (1,0) }$ and ${ (s_1,s_2) =(0,1) }$ so that
\begin{enumerate}
\item[(g)] ${ \bu[y,x_0,x_1,\ldots,x_n] = x_0 \cdot y x_1x_0 \cdot x_2x_1 \cdot x_3x_2 \cdots x_{n-1}x_{n-2} \cdot x_nx_{n-1} y \cdot x_n }$.
\end{enumerate}

Using~(c) instead of~(b), the preceding argument can be repeated to give
\begin{enumerate}
\item[(h)] ${ \bu[z,x_0,x_1,\ldots,x_n] = x_0 \cdot z x_1x_0 \cdot x_2x_1 \cdot x_3x_2 \cdots x_{n-1}x_{n-2} \cdot x_nx_{n-1} z \cdot x_n }$.
\end{enumerate}
It is then easily deduced from~(d), (g), and~(h) that ${ \bu = \bw_n }$.
\end{proof}

\begin{remark}
The requirement that $n \geq 3$ in Lemma~\ref{L: isoterms for S(xyxy)} is necessary since the word~$\bw_2$ is not an isoterm for the monoid $\RQ\{xyxy\}$.
Specifically, $\RQ\{xyxy\}$ satisfies the nontrivial identity ${ \bw_2\approx x_0x_1 \cdot y z \cdot x_0x_2 \cdot yz \cdot x_1x_2 }$.
\end{remark}

\begin{lemma} \label{L: wk xyxy}
Let ${ n,k \geq 3 }$ with ${ n \neq k }$ and let ${ \theta: \sA \to \sA^* }$ be any substitution\up.
\begin{enumerate}
\item[\rm(i)] If ${ \bw_n \theta \preceq \bw_k }$\up, then ${ \bw_n\theta = \bw_n'\theta }$\up.
\item[\rm(ii)] If ${ \bw_n' \theta \preceq \bw_k }$\up, then ${ \bw_n\theta = \bw_n'\theta }$\up.
\end{enumerate}
\end{lemma}

\begin{proof}
Since~$\bw_n$ and~$\bw_n'$ are obtained from one another by in\-ter\-chang\-ing the variables~$y$ and~$z$, it suffices to verify part~(i).
Suppose that ${ \bw_n \theta \preceq \bw_k }$.
Then clearly ${ \bw_n \theta = \bw_n' \theta }$ if ${ y \theta = 1 }$ or ${ z \theta = 1 }$.
Thus assume that ${ y \theta \neq 1 \neq z \theta }$; in what follows, a contradiction is deduced from this assumption, whence the proof is complete.

First observe that
\begin{enumerate}
\item[(a)] the only factors of~$\bw_k$ that have more than one occurrence in~$\bw_k$ are the individual variables and the word~$yz$.
\end{enumerate}
Therefore since ${ \bw_n \theta = \cdots (y \theta)(z \theta) \cdots (y \theta)(z \theta) \cdots \preceq \bw_k }$ with ${ y \theta \neq 1 \neq z \theta }$, the only possibility is ${ y\theta =y }$ and ${ z\theta=z }$.
\[
\begin{picture}(280,50)(00,-05) \setlength{\unitlength}{0.4mm}
\put(06,17){\makebox(0,0)[r]{$\theta:$}}
\put(10,31){\makebox(0,0)[l]{$\bw_n = x_0 \cdot y \ z \cdot x_1x_0 \cdot x_2x_1 \cdots x_nx_{n-1} \cdot y \ z \cdot x_n$}}
\put(10,00){\makebox(0,0)[l]{$\bw_k = x_0 \cdot y \ z \cdot x_1x_0 \cdot x_2x_1 \cdots x_kx_{k-1} \cdot y \ z \cdot x_k$}}
\put(51,26){\vector(0,-1){20}} \put(59,26){\vector(0,-1){20}}
\put(159,26){\vector(0,-1){20}} \put(167,26){\vector(0,-1){20}}
\end{picture}
\]
It follows that the factor ${ \bw_n \theta }$ of~$\bw_k$ contains both occurrences of~$yz$ and so also all variables in between, whence
\begin{enumerate}
\item[(b)] ${ \con(\bw_n\theta) = \con(\bw_k) }$.
\end{enumerate}

Now since ${ y\theta =y }$ and ${ z\theta=z }$, it follows from~(a) that ${ {x_j}{\theta} \in \{1,x_0,x_1,\ldots,x_k\} }$ for all~$x_j$ in~$\con(\bw_n)$.
It thus follows from~(b) that
\begin{enumerate}
\item[(c)] for any ${ x_i \in \con(\bw_k) }$, there exists some ${ x_j \in \con(\bw_n) }$ such that ${ x_j\theta = x_i }$;

\item[(d)] if ${ x_j\theta = x_i }$, then~$\theta$ sends the first (respectively, second)~$x_j$ in~$\bw_n$ to the first (respectively, second)~$x_i$ in~$\bw_k$.
\end{enumerate}
In particular, (c) and the assumption ${ n \neq k }$ imply that ${ n > k }$.
\[
\begin{picture}(345,50)(00,-05) \setlength{\unitlength}{0.4mm}
\put(06,17){\makebox(0,0)[r]{$\theta:$}}
\put(10,30){\makebox(0,0)[l]{$\bw_n = x_0 \cdot y \ z \cdot x_1x_0 \cdot x_2x_1 \cdots x_{k-1}x_{k-2} \cdot x_k \ x_{k-1} \cdot x_{k+1}x_k \cdots x_nx_{n-1} \cdot y \ z \cdot x_n$}}
\put(10,00){\makebox(0,0)[l]{$\bw_k = x_0 \cdot y \ z \cdot x_1x_0 \cdot x_2x_1 \cdots x_{k-1}x_{k-2} \cdot x_k \ x_{k-1} \cdot y \ z \cdot x_k$}}
\put(51,26){\vector(0,-1){20}} \put(59,26){\vector(0,-1){20}}
\put(205,20){\vector(0,-1){14}} \put(205,20){\line(1,0){76}} \put(281,26){\line(0,-1){6}}
\put(212,14){\vector(0,-1){8}} \put(212,14){\line(1,0){76}} \put(288,26){\line(0,-1){12}}
\end{picture}
\]

By~(c), there exists some ${ x_j \in \con(\bw_n) }$ such that ${ {x_j}\theta = x_0 }$.
By~(d), the substitution~$\theta$ sends the first (respectively, second)~$x_j$ in~$\bw_n$ to the first (respectively, second)~$x_0$ in~$\bw_k$.
Since the first factor~$yz$ of~$\bw_k$ is sandwiched between the two occurrences of~$x_0$, the first factor~$yz$ of~$\bw_n$ is also sandwiched between the two occurrences of~$x_j$; this forces ${ j=0 }$, so that ${ x_0 \theta = x_0 }$.
\[
\begin{picture}(345,50)(00,-05) \setlength{\unitlength}{0.4mm}
\put(06,17){\makebox(0,0)[r]{$\theta:$}}
\put(10,30){\makebox(0,0)[l]{$\bw_n = x_0 \cdot y \ z \cdot x_1x_0 \cdot x_2x_1 \cdots x_{k-1}x_{k-2} \cdot x_k \ x_{k-1} \cdot x_{k+1}x_k \cdots x_nx_{n-1} \cdot y \ z \cdot x_n$}}
\put(10,00){\makebox(0,0)[l]{$\bw_k = x_0 \cdot y \ z \cdot x_1x_0 \cdot x_2x_1 \cdots x_{k-1}x_{k-2} \cdot x_k \ x_{k-1} \cdot y \ z \cdot x_k$}}
\put(37,26){\vector(0,-1){20}} \put(51,26){\vector(0,-1){20}} \put(59,26){\vector(0,-1){20}} \put(79,26){\vector(0,-1){20}}
\put(205,20){\vector(0,-1){14}} \put(205,20){\line(1,0){76}} \put(281,26){\line(0,-1){6}}
\put(212,14){\vector(0,-1){8}} \put(212,14){\line(1,0){76}} \put(288,26){\line(0,-1){12}}
\end{picture}
\]

Since ${ k < n }$, the argument in the previous paragraph can be repeated to show that ${ x_i \theta = x_i }$ for all subsequent ${ i = 1,2, \ldots, k-1 }$.
\[
\begin{picture}(345,50)(00,-05) \setlength{\unitlength}{0.4mm}
\put(06,17){\makebox(0,0)[r]{$\theta:$}}
\put(10,30){\makebox(0,0)[l]{$\bw_n = x_0 \cdot y \ z \cdot x_1x_0 \cdot x_2x_1 \cdots x_{k-1}x_{k-2} \cdot x_k \ x_{k-1} \cdot x_{k+1}x_k \cdots x_nx_{n-1} \cdot y \ z \cdot x_n$}}
\put(10,00){\makebox(0,0)[l]{$\bw_k = x_0 \cdot y \ z \cdot x_1x_0 \cdot x_2x_1 \cdots x_{k-1}x_{k-2} \cdot x_k \ x_{k-1} \cdot y \ z \cdot x_k$}}
\put(37,26){\vector(0,-1){20}} \put(51,26){\vector(0,-1){20}} \put(59,26){\vector(0,-1){20}} \put(70,26){\vector(0,-1){20}} \put(79,26){\vector(0,-1){20}}
\put(94,26){\vector(0,-1){20}} \put(103,26){\vector(0,-1){20}} \put(125,26){\vector(0,-1){20}} \put(144,26){\vector(0,-1){20}} \put(181,26){\vector(0,-1){20}}
\put(205,20){\vector(0,-1){14}} \put(205,20){\line(1,0){76}} \put(281,26){\line(0,-1){6}}
\put(212,14){\vector(0,-1){8}} \put(212,14){\line(1,0){76}} \put(288,26){\line(0,-1){12}}
\end{picture}
\]
By~(c), there exists some ${ x_j \in \con(\bw_n) }$ such that ${ {x_j}\theta = x_k }$.
Since the factor~$x_{k-1}yz$ of~$\bw_k$ is sandwiched between the two occurrences of~$x_k$, both the second variable~$x_{k-1}$ in~$\bw_n$ and the second factor~$yz$ in~$\bw_n$ are sandwiched between the two occurrences of~$x_j$.
But it is easily seen that no such variable~$x_j$ exists.
\end{proof}

\begin{theorem}\label{T: xyxy}
The variety $\bbM\{xyxy\}$ is of~\two\up.
\end{theorem}

\begin{proof}
Let~$\Gamma$ be any identity basis for the monoid $\RQ\{xyxy\}$.
For any subset~$\mathtt{N}$ of ${ \{ 3,4,5,\ldots\} }$, define the identity system ${ \Gamma_\mathtt{N} = \{ \bw_n \approx \bw_n' \,|\, n \in \mathtt{N} \} }$.
In what follows, it is shown that the deduction ${ \Gamma \cup \Gamma_\mathtt{N} \vdash \bw_k \approx \bw_k' }$ holds if and only if ${ k \in \mathtt{N} }$.
Therefore each subset of ${ \{3,4,5,\ldots\} }$ corresponds to a unique subvariety of $\bbM\{xyxy\}$, whence the theorem is established.

Suppose that the deduction ${ \Gamma \cup \Gamma_\mathtt{N} \vdash \bw_k \approx \bw_k' }$ holds for some ${ k \notin \mathtt{N} }$.
Then there exists a sequence
\[
\bw_k=\bu_0, \bu_1, \ldots, \bu_r = \bw_k'
\]
of distinct words such that each identity ${ \bu_i \approx \bu_{i+1} }$ is directly deducible from some identity~$\gamma_i$ in ${ \Gamma \cup \Gamma_\mathtt{N} }$.
By Lemma~\ref{L: isoterms for S(xyxy)}, the word ${ \bw_k = \bu_0 }$ is an isoterm for $\RQ\{xyxy\}$.
Therefore the monoid $\RQ\{xyxy\}$ does not satisfy the identity ${ \bu_0 \approx \bu_1 }$, whence ${ \Gamma \nvdash \bu_0 \approx \bu_1 }$.
Hence the identity~$\gamma_0$ is in~$\Gamma_\mathtt{N}$, say~$\gamma_0$ is ${ \bw_n \approx \bw_n' }$ for some ${ n \in \mathtt{N} }$.
Now ${ \{\bu_0,\bu_1\}=\{\ba(\bw_n\theta_0)\bb,\ba(\bw_n'\theta)\bb\} }$ for some ${ \ba,\bb\in\sA^* }$ implies that either $\bw_n\theta_0$ or $\bw_n'\theta_0$ is a factor of ${ \bu_0 = \bw_k }$.
Since ${ n \neq k }$, it follows from Lemma~\ref{L: wk xyxy} that ${ \bw_n \theta_0 = \bw_n' \theta_0 }$, whence the contradiction ${ \bu_0 = \bu_1 }$ is deduced.
\end{proof}

As shown in Jackson and~Sapir~\cite{JS00}, there exist infinitely many finite sets~$\sW$ of words such that $\RQ\sW$ is finitely based with ${ \RQ\{xyxy\} \in \bbM\sW }$; one such example is ${ \sW=\{xyxy,xy^2x,x^2y^2\} }$.
Therefore there exists an abundance of finitely based, finitely generated monoid varieties of \two.

Many finite $\mathscr{J}$-trivial monoids from Volkov~\cite{Vol04}, several of which are finitely based, also generate monoid varieties of {\two} due to the inclusion of $\RQ\{xyxy\}$.

\subsection{Joins of Cross monoid varieties} \label{subsec: join}

\begin{lemma}[Oates and Powell~\cite{OP64}] \label{L: groups}
Any finite group generates a Cross monoid variety\up.
\end{lemma}

\begin{lemma}[Lee~{\cite[Lemma~3.4(i) and Proposition~4.1]{Lee12}}] \label{L: Lee12}
Any monoid that satisfies the identity ${ xyxzx \approx xyzx }$ generates a Cross monoid variety\up.
\end{lemma}

\begin{lemma}[Lee and Zhang~\cite{LZ14}] \label{L: min non small}
Up to isomorphism and anti-isomorphism\up, the monoid~$P_2^1$ obtained from the semigroup
\[
P_2 = \big\langle\, a,b \,\big|\, a^2=ab=a, \, b^2a=b^2 \,\big\rangle
\]
is the only monoid of order five that generates a non\-small monoid variety\up; every other monoid of order five or less generates a Cross monoid variety\up.
\end{lemma}

Using Theorem~\ref{T: xyxy} and Lemmas~\ref{L: groups}--\ref{L: min non small}, many examples of two Cross monoid varieties can be found with join that is of~\two.
For instance, the symmetric group~$S_3$ on three symbols and the monoid~$\RQ\{x^2\}$ of order four generate Cross monoid varieties.
But since the word~$xyxy$ is an isoterm for ${ \big\{S_3, \RQ\{x^2\}\big\} }$, the join
\[
\Vm\{S_3\} \vee \bbM\{x^2\} = \Vm\big\{S_3, \RQ\{x^2\}\big\}
\]
contains $\RQ\{xyxy\}$ by Lemma~\ref{L: isoterm} and so is of~\two.
For an aperiodic example, consider the semigroup
\[
A_0 = \big\langle\, a,b \,\big|\, a^2=a, \, b^2=b, \, ab=0\big\rangle
\]
of order four.
Then the variety ${ \Vm\{L_2^1,A_0^1,R_2^1\} }$ is Cross by Lemma~\ref{L: Lee12}, but the word~$xyxy$ is an isoterm for ${ \big\{L_2^1,A_0^1,R_2^1, \RQ\{x^2\}\big\} }$ so that the join
\[
\Vm\{L_2^1,A_0^1,R_2^1\} \vee \bbM\{x^2\} = \Vm\big\{L_2^1,A_0^1,R_2^1, \RQ\{x^2\}\big\}
\]
contains $\RQ\{xyxy\}$ by Lemma~\ref{L: isoterm}.

Presently, a few questions regarding semigroup varieties remain open.

\begin{question}[Jackson~{\cite[Question~3.15]{Jac00}}] \label{Q: Jackson}
Are there finitely generated, small semigroup varieties~$\cU$ and~$\cV$ such that the join ${ \cU \vee \cV }$ is of~\two?
\end{question}

A variety~$\cV$ is said to satisfy the \textit{ascending chain condition} if the lattice $\fL(\cV)$ satisfies the same condition.

\begin{question}[Shevrin et al.~{\cite[Question~10.2(a)]{SVV09}}]  \label{Q: SVV09}
Are there semigroup varieties~$\cU$ and~$\cV$ with the ascending chain condition such that the join ${ \cU \vee \cV }$ does not satisfy the same condition?
\end{question}

Observe that the proof of Theorem~\ref{T: xyxy} actually established a stronger result: the lattice of all subsets of ${ \{3,4,5,\ldots\} }$ order-embeds into the lattice $\fL(\bbM\{xyxy\})$.
Hence within the context of monoids, Questions~\ref{Q: Jackson} and~\ref{Q: SVV09} are positively answered by the Cross varieties ${ \cU \in \big\{ \Vm\{L_2^1,A_0^1,R_2^1\}, \Vm\{S_3\} \big\} }$ and ${ \cV = \bbM\{x^2\} }$.

\subsection{Variety generated by all monoids of order~$n$} \label{subsec: order n}

Let~$\cM_n$ denote the class of all monoids of order~$n$.
Since there exist finitely many pairwise nonisomorphic monoids of order~$n$, the variety $\Vm\cM_n$ is the join of finitely many varieties.
By Lemma~\ref{L: min non small}, the variety $\Vm\cM_4$ is the join of finitely many Cross varieties.

\begin{proposition} \quad
\begin{enumerate}
\item[\rm(i)] The variety $\Vm\cM_2$ contains four subvarieties\up.
\item[\rm(ii)] The variety $\Vm\cM_3$ contains~$60$ subvarieties\up.
\item[\rm(iii)] The variety $\Vm\cM_4$ is of~\two\up.
\end{enumerate}
\end{proposition}

\begin{proof}
For each ${ n \geq 1 }$, let~$Z_n$ denote the cyclic group of order~$n$.

(i) It is well known that the variety ${ \Vm\cM_2 = \Vm\{\RQ\varnothing,Z_2\} }$ contains four subvarieties:~$\mathbf{0}$, $\bbM\varnothing$, $\Vm\{Z_2\}$, and itself.

(ii) Up to isomorphism and anti-isomorphism, all semigroups of order three are listed in Luo and Zhang~\cite[Table~2]{LuoZ11}.
It is easily deduced that \[\Vm\cM_3 = \Vm\{L_2^1, \RQ\{x\}, R_2^1, Z_2, Z_3\}.\]
Following Luo and Zhang~\cite[\S4]{LuoZ11}, the variety ${ \Vm\{L_2^1, \RQ\{x\}, R_2^1, Z_n\} }$ can be shown to be defined by the identities
\[
x^{n+1}hx \approx xhx, \quad xhxtx \approx x^2htx, \quad xhxyty \approx xhyxty.
\]
It is then routinely shown that
\begin{enumerate}
\item[(a)] the lattice $\fL(\Vm\cM_3)$ is the disjoint union of $\fL(\Vm\{L_2^1, \RQ\{x\}, R_2^1\})$,
\begin{align*}
& [\Vm\{Z_2\}, \Vm\{L_2^1, \RQ\{x\}, R_2^1,Z_2\}], \\
& [\Vm\{Z_3\}, \Vm\{L_2^1, \RQ\{x\}, R_2^1,Z_3\}], \\
\text{and} \quad & [\Vm\{Z_2,Z_3\}, \Vm\{L_2^1, \RQ\{x\}, R_2^1,Z_2,Z_3\}];
\end{align*}

\item[(b)] the latter three intervals in~(a) are isomorphic to ${ \fL(\Vm\{L_2^1, \RQ\{x\}, R_2^1\}) }$.
\end{enumerate}
By Proposition~\ref{P: LNR}, the lattice ${ \fL(\Vm\{L_2^1, \RQ\{x\}, R_2^1\}) }$ contains~$15$ varieties.
Therefore by~(a) and~(b), the number of subvarieties of $\Vm\cM_3$ is ${ 15 \times 4 = 60 }$.

(iii) The class~$\cM_4$ contains~$L_2^1$, $\RQ\{x^2\}$, and the monoid~$T_2$ of transformations of ${ \{ 1,2 \} }$, which can be given by
\[
T_2 = \big\langle\, a,b \,\big|\, a^2=ba=a, \, b^2=1 \,\big\rangle.
\]
Hence ${ \Vm\{ L_2^1, \RQ\{x^2\},T_2\} \subseteq \Vm\cM_4 }$.
It is easily checked that~$xyxy$ is an isoterm for ${ \{ L_2^1, \RQ\{x^2\},T_2 \} }$ so that ${ \RQ\{xyxy\} \in \Vm\{ L_2^1, \RQ\{x^2\},T_2 \} }$ by Lemma~\ref{L: isoterm}.
The result thus holds by Theorem~\ref{T: xyxy}.

Alternately, it follows from Lee and Li~\cite[Theorem~1.1]{LL15} that the identities
\begin{gather*} \label{basis M4}
x^{13} hxkx \approx xhxkx, \quad xhx^2kx \approx x^3hkx, \quad xhy^2x^2ky \approx xhx^2y^2ky, \\
xhykxytxdy \approx xhykyxtxdy, \quad xhykxytydx \approx xhykyxtydx
\end{gather*}
constitute a basis for the variety $\Vm\cM_4$; the monoid $\RQ\{xyxy\}$ satisfies these identities and so belongs to the variety.
\end{proof}

Consequently, the monoid variety $\Vm\cM_n$ is of~{\two} if and only if ${ n \geq 4 }$.
It turns out that by Jackson~\cite[Corollary~3.17]{Jac00} and Luo and Zhang~\cite{LuoZ11}, a similar result holds for the semigroup variety $\Vs\cM_n$.
\begin{gather*}
\begin{tabular}{rccccccc} \hline
                                     & &       & &            & &                \\[-0.1in]
                                     & & $n=2$ & & $n=3$      & & $n\geq4$       \\[-0.1in]
                                     & &       & &            & &                \\ \hline\hline
                                     & &       & &            & &                \\[-0.1in]
Number of subvarieties of $\Vs\cM_n$ & & $32$  & & $\aleph_0$ & & $2^{\aleph_0}$ \\
                                     & &       & &            & &                \\[-0.1in]
Number of subvarieties of $\Vm\cM_n$ & & $4$   & & $60$       & & $2^{\aleph_0}$ \\
                                     & &       & &            & &                \\[-0.1in] \hline
\end{tabular}
\end{gather*}

\section{Varieties generated by {\infb} finite monoids} \label{sec: INFB}

An algebra~$A$ is \textit{\infb} if every lo\-cally finite variety containing~$A$ is {\nfb}.
M.\,V.~Sapir~\cite{mS88b} provided an elegant description of {\infb} finite semigroups based on \textit{Zimin words} ${ \bz_1,\bz_2,\bz_3,\ldots}$ defined by ${ \bz_1 = x_1 }$ and ${ \bz_{n+1} =\bz_n x_{n+1}\bz_n }$ for all ${ n\geq 1 }$.

\begin{theorem}[M.\,V.~Sapir~{\cite[Proposition~7]{mS88b}}] \label{T: Sapir}
A finite semigroup is {\infb} if and only if all Zimin words are isoterms for it\up.
\end{theorem}

Using Theorem~\ref{T: Sapir}, Jackson \cite[proof of Corollary~3.9]{Jac00} demonstrated that the variety generated by any {\infb} finite semigroup contains the monoid $\RQ\{xyx\}$; in view of Proposition~\ref{P: xyx}, such a variety is of~\two.
But since the monoid variety $\bbM\{xyx\}$ contains only five subvarieties~\cite[Lemma~4.4]{Jac05b}, the inclusion of the monoid $\RQ\{xyx\}$ is insufficient for a monoid variety to be of~\two.

M.\,V.~Sapir has also provided a structural description of {\infb} finite semigroups that is detailed enough to provide an algorithm that decides if a finite semigroup is {\infb} \cite[Theorem~1]{mS88a}.
It is natural to question if such an algorithm exists for finite monoids.

\begin{question}[M.\,V.~Sapir~{\cite[Problem~3.10.12]{mS14}}] \label{Q: Sapir}
Is there an algorithm that decides if a finite monoid is {\infb} as a monoid?
\end{question}

In the present section, arguments from~\S\ref{sec: xyxy} are adapted to show that any {\infb} finite monoid generates a monoid variety of~\two.
In~\S\ref{subsec: description}, it is shown that a finite monoid is {\infb} as a semigroup if and only if it is {\infb} as a monoid.
It follows that Question~\ref{Q: Sapir} is affirmatively answered by the aforementioned algorithm of M.\,V.~Sapir \cite[Theorem~1]{mS88a}.
In~\S\ref{subsec: Zimin}, some preliminary results regarding Zimin words are established.
Main results of the section are then established in~\S\ref{subsec: infb}.

\subsection{Inherently {\nfb} finite monoids} \label{subsec: description}

\begin{theorem} \label{T: infb equivalence}
The following statements on any finite monoid~$M$ are equivalent\up:
\begin{enumerate}
\item[\rm(i)] $M$ belongs to some finitely based locally finite semigroup variety\up;
\item[\rm(ii)] $M$ belongs to some finitely based locally finite monoid variety\up.
\end{enumerate}
Consequently\up, $M$ is an {\infb} semigroup if and only if~$M$ is an {\infb} monoid\up.
\end{theorem}

\begin{proof}
(i)~implies~(ii).
Suppose that~$M$ belongs to some finitely based locally finite semigroup variety~$\bV$.
Let~$\Sigma$ be any finite identity basis for~$\bV$ and let~$\bbV$ denote the monoid variety defined by~$\Sigma$.
Evidently, the monoid variety~$\bbV$ is finitely based and ${ M \in \bbV }$.
Moreover, $\bbV$ is locally finite, because any of its finitely gen\-er\-ated monoids has its semigroup reduct in~$\bV$, and hence is finite.

(ii)~implies~(i).
Suppose that~$M$ belongs to some finitely based locally finite monoid variety~$\bbV$.
Let~$\Sigma$ be any finite identity basis for~$\bbV$.
If~$\Sigma$ contains some identity ${ \bu \approx \bv }$ with ${ \con(\bu) \neq \con(\bv) }$, say ${ x \in \con(\bu) \backslash \con(\bv) }$, then~$M$ satisfies the identity ${ x^{\occ(x,\bu)} \approx 1 }$ and so is a finite group; it is well known that all finite groups are finitely based~\cite{OP64}, so that $\Vs\{M\}$ is the required finitely based locally finite semigroup variety containing~$M$.
Therefore it suffices to assume that every identity ${ \bu \approx \bv }$ in~$\Sigma$ satisfies the property that ${ \con(\bu) = \con(\bv) }$.
Let~$\Sigma'$ be the result of adding to~$\Sigma$ all nontrivial identities obtained from identities in~$\Sigma$ by assigning~$1$ to some variables and reducing modulo ${ x1\mapsto x }$ and ${ 1x\mapsto x }$.
It is clear that in the semigroup signature, $\Sigma'$ has the same power of deduction for semigroup identities as~$\Sigma$ has for semigroup identities in the monoid signature: any application of $\bu\approx \bv$ from~$\Sigma$ in an equational deduction can be replaced by an application of some identity in~$\Sigma'$ obtained from ${ \bu\approx \bv }$ by variable deletion.
Thus~$\Sigma'$ is a finite identity basis defining the semigroup variety~$\bV$ gen\-er\-ated by the semigroup reducts of monoids in~$\bbV$.
Clearly the semigroup variety $\bV$ is finitely based and ${ M \in \bV }$.
As the monoid variety~$\bbV$ is locally finite, it is also \textit{uniformly locally finite} in the sense that for each ${ n \in \{1,2,3,\ldots \} }$, there is a number~$k$ such that any $n$-gen\-er\-ated member has at most~$k$ elements~\cite{Mal73}.
Therefore the semigroup variety~$\bV$ is also locally finite.
\end{proof}

\subsection{Some results on Zimin words} \label{subsec: Zimin}

\begin{lemma} \label{L: Zn+1 isoterm}
Suppose that~$\bz_{n+1}$ is an isoterm for a monoid~$M$\up.
Then~$\bz_n$ is also an isoterm for~$M$\up.
\end{lemma}

\begin{proof}
If~$\bz_n$ is not an isoterm for~$M$, say~$M$ satisfies a nontrivial identity ${ \bz_n \approx \bw }$ for some ${ \bw \in \sA^+ }$, then~$M$ also satisfies the nontrivial identity ${ \bz_n x_{n+1} \bz_n \approx \bw x_{n+1} \bw }$, whence ${ \bz_{n+1} = \bz_n x_{n+1} \bz_n }$ is not an isoterm for~$M$.
\end{proof}

\begin{lemma} \label{L: Z2 isoterm}
Let ${ \bw \approx \bw' }$ be any identity satisfied by a monoid~$M$\up.
Suppose that~$\bz_n$ is an isoterm for~$M$ for some ${ n \geq 2 }$\up.
Then\up:
\begin{enumerate}
\item[\rm(i)] ${ \con(\bw) = \con(\bw') }$\up;
\item[\rm(ii)] ${ \occ(x,\bw) = 1 }$ if and only if ${ \occ(x,\bw') = 1 }$\up;
\item[\rm(iii)] ${ \bw[x,y] = xyx }$ if and only if ${ \bw'[x,y] = xyx }$\up.
\end{enumerate}
\end{lemma}

\begin{proof}
By Lemma~\ref{L: Zn+1 isoterm}, the words ${ \bz_1 = x_1 }$ and ${ \bz_2 = x_1x_2x_1 }$ are isoterms for~$M$.
Parts~(i) and~(ii) are then easily established.
Part~(iii) holds because the identity ${ \bw[x,y] \approx \bw'[x,y] }$ is satisfied by~$M$ and~$xyx$ is an isoterm.
\end{proof}

\begin{lemma} \label{L: zn form}
For each ${ n\geq3 }$\up, the word~$\bz_n$ can be written as
\begin{equation}
\bz_n = \bp_1 \bigg(\prod_{i=1}^{n-1}(\bp_{i+1} \bp_i)\bigg) \bq_n = \bp_1 \cdot \bp_2 \bp_1 \cdots \bp_{n-1} \bp_{n-2} \cdot \bp_n \bp_{n-1} \cdot \bq_n, \label{Zn}
\end{equation}
where ${ \bp_1,\bp_2,\ldots,\bp_n,\bq_n \in \sA^+ }$ satisfy the following\up:
\begin{enumerate}
\item[\rm(I)] ${ \con(\bp_i) \subseteq \{ x_1,x_2,\ldots,x_i\} }$ for each ${ i \in \{1,2,\ldots,n \} }$\up;
\item[\rm(II)] ${ \occ(x_i,\bp_i) = 1 }$ for each ${ i \in \{1,2,\ldots,n \} }$\up;
\item[\rm(III)] ${ \con(\bq_n) \subseteq \{x_1,x_2,\ldots,x_{n-2} \} }$\up.
\end{enumerate}
\end{lemma}

\begin{proof}
The basic case holds because
\[
\bz_3 = \underbrace{x_1}_{\bp_1} \ \underbrace{x_2}_{\bp_2} \ \underbrace{x_1}_{\bp_1} \ \underbrace{x_3x_1}_{\bp_3} \ \underbrace{x_2}_{\bp_2} \ \underbrace{x_1}_{\bq_3}.
\]
Suppose that ${ n \geq 4 }$ and the result holds for~$\bz_{n-1}$. That is,
\[
\bz_{n-1} =  \bp_1 \cdot \bp_2 \bp_1 \cdots \bp_{n-1} \bp_{n-2} \cdot \bq_{n-1}
\]
where ${ \bp_1,\bp_2,\ldots,\bp_{n-1},\bq_{n-1} \in \sA^+ }$ satisfy the following:
\begin{enumerate}
\item[\rm(I$'$)] ${ \con(\bp_i) \subseteq \{ x_1,x_2,\ldots,x_i \} }$ for each ${ i \in \{1,2,\ldots,n-1 \} }$;
\item[\rm(II$'$)] ${ \occ(x_i,\bp_i)=1 }$ for each ${ i \in \{1,2,\ldots,n-1 \} }$;
\item[\rm(III$'$)] ${ \con(\bq_{n-1}) \subseteq \{x_1,x_2,\ldots,x_{n-3} \} }$.
\end{enumerate}
Then
\begin{align*}
\bz_n & = \bz_{n-1} \, x_n \, \bz_{n-1} \\
& = \bp_1 \cdot \bp_2 \bp_1 \cdots \bp_{n-1} \bp_{n-2} \cdot \underbrace{\bq_{n-1} \, x_n \, \bp_1 \cdot \bp_2 \bp_1 \cdots \bp_{n-2} \bp_{n-3}}_{\bp_n} \,\cdot \, \bp_{n-1} \underbrace{\bp_{n-2} \cdot \bq_{n-1}}_{\bq_n}
\end{align*}
is of the form~\eqref{Zn}. It follows from (I$'$)--(III$'$) that
\begin{align*}
\con(\bp_n) & = \con(\bq_{n-1}) \cup \{x_n\} \cup \con(\bp_1 \cdot \bp_2 \bp_1 \cdots \bp_{n-2} \bp_{n-3}) \\
& \subseteq \{x_1,x_2,\ldots,x_{n-3}\} \cup \{x_n\} \cup \{ x_1,x_2,\ldots,x_{n-2}\} \\
& \subseteq \{x_1,x_2,\ldots,x_n\}, \\
\occ(x_n,\bp_n) & = \occ(x_n,\bq_{n-1})+1 + \occ(x_n,\bp_1 \cdot \bp_2 \bp_1 \cdots \bp_{n-2} \bp_{n-3}) \\
& = 0+1+0 = 1, \\
\text{and} \quad \con(\bq_n) & = \con(\bp_{n-2}) \cup \con(\bq_{n-1}) \\
& \subseteq \{ x_1,x_2,\ldots,x_{n-2}\} \cup \{ x_1,x_2,\ldots,x_{n-3}\} \\
& \subseteq \{ x_1,x_2,\ldots,x_{n-2}\}.
\end{align*}
Therefore~(I)--(III) hold.
\end{proof}

\subsection{Main results of section} \label{subsec: infb}

\begin{theorem} \label{T: all Zimin isoterms}
Suppose that all Zimin words are isoterms for a monoid $M$\up. Then the monoid variety $\Vm\{M\}$ is of~\two\up.
\end{theorem}

\begin{corollary} \label{C: var Zimin}
The monoid variety ${ \bbM\{\bz_1, \bz_2,\ldots\} }$ is of~\two\up.
\end{corollary}

\begin{proof}
By Lemma~\ref{L: isoterm}, all Zimin words are isoterms for ${ \bbM\{\bz_1, \bz_2,\ldots\} }$.
The result then follows from Theorem~\ref{T: all Zimin isoterms}.
\end{proof}

\begin{corollary} \label{C: infb}
Any {\infb} finite monoid generates a monoid variety of~\two\up.
\end{corollary}

\begin{proof}
By Theorems~\ref{T: Sapir} and~\ref{T: infb equivalence}, all Zimin words are isoterms for any {\infb} finite monoid.
The result then follows from Theorem~\ref{T: all Zimin isoterms}.
\end{proof}

It is well known that the monoid~$A_2^1$ obtained from
\[
A_2 = \big\langle\, a,b \,\big|\, a^2=aba=a, \, b^2=0, \, bab=b \,\big\rangle
\]
and the monoid~$B_2^1$ are {\infb} semigroups~\cite{mS88b} of order six.
Trahtman~\cite{Tra88} proved that the semigroup variety $\Vs\{A_2^1\}$ is of~\two, while Jackson~\cite{Jac00} proved that the semigroup variety $\Vs\{B_2^1\}$ also has the same type.
(Three finitely based monoids of order six have also been found to generate semigroup varieties of~\two~\cite{ELL10}.)
The question of whether or not the monoid varieties $\Vm\{A_2^1\}$ and $\Vm\{B_2^1\}$ are of~{\two} was posed by Lee~\cite[Question~7.4]{Lee12}.
By Theorem~\ref{T: infb equivalence} and Corollary~\ref{C: infb}, the answer to this question is affirmative.

The remainder of this subsection is devoted to the proof of Theorem~\ref{T: all Zimin isoterms}.
For each ${ n \geq 3 }$, define the words
\begin{align*}
\bw_n & = x_0 h x_1yz x_0 \cdot x_2x_1 \cdot x_3x_2 \cdots x_{n-1}x_{n-2} \cdot x_n yz x_{n-1} t x_n \\
\text{and} \quad \bw_n' & = x_0 h x_1zy x_0 \cdot x_2x_1 \cdot x_3x_2 \cdots x_{n-1}x_{n-2} \cdot x_n zy x_{n-1} t x_n.
\end{align*}

\begin{lemma}\label{L: isoterm 0hyz1}
Suppose that~$\bz_5$ is an isoterm for a monoid~$M$\up. Then
\[
\bw_n [x_0,h,y,z,x_1] \quad \text{and} \quad \bw_n [x_n,t,y,z,x_{n-1}]
\]
are also isoterms for~$M$\up.
\end{lemma}

\begin{proof}
By symmetry, it suffices to show that
\[
\bw = \bw_n [x_0,h,y,z,x_1] = x_0 h x_1 yz x_0x_1yz
\]
is an isoterm for the monoid~$M$.
Let ${ \bw \approx \bw' }$ be any identity satisfied by~$M$.
Since ${ \con(\bw) = \{ x_0,h,y,z,x_1 \} }$ and ${ \bw[x_0,h] = x_0 h x_0 }$, it follows from Lemma~\ref{L: Z2 isoterm} that ${ \con(\bw') = \{ x_0,h,y,z,x_1 \} }$ and ${ \bw'[x_0,h] = x_0 h x_0 }$.
Therefore
\begin{enumerate}
\item[(a)] ${ \bw' = \ba x_0 \bb h \bc x_0 \bd }$ for some ${ \ba,\bb,\bc,\bd \in \{ y,z,x_1\}^* }$.
\end{enumerate}
Let ${ \theta : \{x_0,h,y,z,x_1\} \to \sA^+ }$ denote the substitution given by
\begin{gather*}
(x_0,h,y,z,x_1) \mapsto (x_1x_2x_1x_4x_1, \, x_2x_1x_3x_1x_2x_1x_5x_1, \, x_1, \, x_3, \, x_2).
\end{gather*}
Then
\begin{enumerate}
\item[(b)] ${ \bw\theta = {\underbrace{x_1x_2x_1x_4x_1}_{x_0 \theta}} \cdot {\underbrace{x_2x_1x_3x_1x_2x_1x_5x_1}_{h\theta}} \cdot {\underbrace{x_2 x_1 x_3}_{(x_1 yz)\theta} \cdot \underbrace{x_1x_2x_1x_4x_1}_{x_0 \theta}} \cdot {\underbrace{x_2 x_1 x_3}_{(x_1 yz)\theta}} }$
\end{enumerate}
is a factor of~$\bz_5$, that is, ${ \bz_5 = \mathbf{e} (\bw \theta) \mathbf{f} }$ for some ${ \mathbf{e}, \mathbf{f} \in \sA^* }$.
If ${ \bw \theta \neq \bw' \theta }$, then ${ \mathbf{e} (\bw\theta) \mathbf{f} \neq \mathbf{e} (\bw' \theta) \mathbf{f} }$ so that ${ \bz_5 \approx \mathbf{e} (\bw' \theta) \mathbf{f} }$ is contradictorily a nontrivial identity satisfied by~$M$.
Therefore ${ \bw \theta = \bw' \theta }$, whence by~(a) and~(b),
\begin{align*}
\bw'\theta & = \ba\theta \cdot x_1x_2x_1x_4x_1 \cdot \bb\theta \cdot x_2x_1x_3x_1x_2x_1x_5x_1 \cdot \bc\theta \cdot x_1x_2x_1x_4x_1 \cdot \bd\theta \\
& = x_1x_2x_1x_4x_1 \cdot x_2x_1x_3x_1x_2x_1x_5x_1 \cdot x_2 x_1 x_3 \cdot x_1x_2x_1x_4x_1 \cdot x_2 x_1 x_3.
\end{align*}
Since ${ \ba,\bb,\bc,\bd \in \{ y,z,x_1\}^* }$ implies that ${ \ba \theta, \bb \theta, \bc \theta, \bd \theta \in \{ x_1,x_2,x_3\}^* }$, it follows that ${ \ba = \bb = 1 }$ and ${ \bc = \bd = x_1yz }$.
Consequently, ${ \bw' \stackrel{\rm(a)}{=} x_0 h x_1 yz x_0x_1yz = \bw }$.
\end{proof}

\begin{lemma} \label{L: isoterm h0...n}
Suppose that~$\bz_{n+2}$ is an isoterm for a monoid~$M$\up.
Then
\[
\bw_n [h,x_0,x_1,\ldots,x_n] \quad \text{and} \quad \bw_n [x_0,x_1,\ldots,x_n,t]
\]
are also isoterms for~$M$\up.
\end{lemma}

\begin{proof}
By symmetry, it suffices to show that
\[
\bw = \bw_n [x_0,x_1,\ldots,x_n,t] = x_0 \cdot x_1 x_0 \cdot x_2x_1 \cdots x_n x_{n-1} \cdot tx_n
\]
is an isoterm for the monoid~$M$.
Let ${ \bw \approx \bw' }$ be any identity satisfied by~$M$.
Since ${ \con(\bw) = \{ x_0,x_1,\ldots,x_n,t \} }$ and ${ \bw[x_n,t] = x_n t x_n }$, it follows from Lemma~\ref{L: Z2 isoterm} that ${ \con(\bw') = \{ x_0,x_1,\ldots,x_n,t \} }$ and ${ \bw'[x_n,t] = x_n t x_n }$.
Therefore
\begin{enumerate}
\item[(a)] ${ \bw' = \ba x_n \bb t \bc x_n \bd }$ for some ${ \ba,\bb,\bc,\bd \in \{ x_0,x_1,\ldots,x_{n-1}\}^* }$.
\end{enumerate}
On the other hand, by Lemma~\ref{L: zn form},
\begin{equation*}
\bz_{n+2} = \bp_1 \cdot \bp_2 \bp_1 \cdot \bp_3 \bp_2 \cdots \bp_{n+1} \bp_n \cdot \bp_{n+2} \bp_{n+1} \cdot \bq_{n+2}
\end{equation*}
for some ${ \bp_1,\bp_2,\ldots,\bp_{n+2},\bq_{n+2} \in \sA^+ }$ such that
\begin{enumerate}
\item[(b)] ${ \con(\bp_i) \subseteq \{ x_1,x_2,\ldots,x_i \} }$ for each ${ i \in \{1,2,\ldots,n+2 \} }$,

\item[(c)] ${ \occ(x_i,\bp_i) = 1 }$ for each ${ i \in \{1,2,\ldots,n+2 \} }$, and

\item[(d)] ${ \con(\bq_{n+2}) \subseteq \{x_1,x_2,\ldots,x_n \} }$.
\end{enumerate}

Let ${ \theta: \{x_0,x_1,\ldots,x_n,t\} \to \sA^+ }$ denote the substitution given by
\[
(x_0,x_1,\ldots,x_{n-1},x_n,t) \mapsto (\bp_1, \bp_2,\ldots, \bp_n, \bp_{n+1}, \bp_{n+2}).
\]
Then
\begin{enumerate}
\item[(e)] ${ \bw \theta = \bp_1 \cdot \bp_2 \bp_1 \cdot \bp_3 \bp_2 \cdots \bp_n \bp_{n-1} \cdot \bp_{n+1} \bp_n \cdot \bp_{n+2} \bp_{n+1} }$ and

\item[(f)] ${ \bw'\theta = \ba \theta \cdot \bp_{n+1} \cdot \bb \theta \cdot \bp_{n+2} \cdot \bc \theta \cdot \bp_{n+1} \cdot \bd \theta }$ where

\item[(g)] ${ \ba\theta, \bb\theta, \bc\theta, \bd\theta \in \{ \bp_1,\bp_2,\ldots,\bp_n\}^* }$.
\end{enumerate}
Note that ${ \bz_{n+2} = (\bw \theta)\bq_{n+2} }$.
Hence if ${ \bw \theta \neq \bw' \theta }$, then ${ (\bw \theta)\bq_{n+2} \neq (\bw' \theta)\bq_{n+2} }$ so that ${ \bz_{n+2} \approx (\bw' \theta)\bq_{n+2} }$ is contradictorily a nontrivial identity satisfied by~$M$.
Therefore ${ \bw \theta = \bw' \theta }$.
Since ${ x_{n+1},x_{n+2} \notin \con\big((\ba\bb\bc\bd)\theta\big) }$ by~(b) and~(g), a simple inspection of~(e) and~(f) yields
\[
\ba \theta = \bp_1 \cdot \bp_2 \bp_1 \cdot \bp_3 \bp_2 \cdots \bp_{n-1} \bp_{n-2} \cdot \bp_n \bp_{n-1}, \ \bb \theta = \bp_n, \ \text{and} \ \bc\theta = \bd\theta = 1.
\]
It is then clear that
\begin{enumerate}
\item[(h)] ${ \bc = \bd = 1 }$ and ${ \bb = x_{n-1} }$.
\end{enumerate}
Since ${ x_n \in \con(\bp_n) \backslash \con(\bp_1 \bp_2 \cdots \bp_{n-1}) }$ by~(b) and~(c), the variable~$x_n$ occurs precisely once in~$\ba\theta$, whence
\[
\ba = \ba' x_{n-1} \bb'
\]
for some ${ \ba',\bb' \in \{x_0,x_1,\ldots,x_{n-2}\}^* }$.
Since the factor~$\bp_n$ of~$\ba\theta$ is preceded by one~$\bp_{n-1}$ and immediately followed by one~$\bp_{n-1}$, it follows from~(b) and~(c) that
\[
\ba = \ba'' x_{n-2} \bb'' x_{n-1} x_{n-2}
\]
for some ${ \ba'',\bb'' \in \{x_0,x_1,\ldots,x_{n-3}\}^* }$.
This argument can be repeated to obtain
\[
\ba = x_0 \cdot x_1 x_0 \cdot x_2x_1 \cdots x_{n-3} x_{n-4} \cdot x_{n-2} x_{n-3} \cdot x_{n-1} x_{n-2}.
\]
Hence ${ \bw' = \ba x_n \bb t \bc x_n \bd = \bw }$ by~(a) and~(h).
\end{proof}

\begin{lemma} \label{L: isoterms for S(zn)}
Suppose that all Zimin words are isoterms for a monoid~$M$\up.
Then for each ${ n \geq 3 }$\up, the words~$\bw_n$ and~$\bw_n'$ are also isoterms for $M$\up.
\end{lemma}

\begin{proof}
By Lemmas~\ref{L: isoterm 0hyz1} and~\ref{L: isoterm h0...n}, all neighboring variables in~$\bw_n$ and~$\bw_n'$ are fixed in their relative positions.
\end{proof}

\begin{lemma} \label{L: wk infb}
Let ${ n,k \geq 3 }$ with ${ n \neq k }$ and let ${ \theta: \sA \to \sA^* }$ be any substitution\up.
\begin{enumerate}
\item[\rm(i)] If ${ \bw_n \theta \preceq \bw_k }$\up, then ${ \bw_n\theta = \bw_n'\theta }$\up.

\item[\rm(ii)] If ${ \bw_n' \theta \preceq \bw_k }$\up, then ${ \bw_n\theta = \bw_n'\theta }$\up.
\end{enumerate}
\end{lemma}

\begin{proof}
This is very similar to Lemma~\ref{L: wk xyxy}.
\end{proof}

\begin{proof}[Proof of Theorem~\ref{T: all Zimin isoterms}]
This is similar to the proof of Theorem~\ref{T: xyxy}.
\end{proof}

\section{A monoid variety of~\one} \label{sec: E}

This section is concerned with the monoid~$E^1$ obtained from the semigroup
\[
E = \big\langle a,b,c \,\big|\, a^2=ab=0, \, ba=ca=a, \, b^2=bc=b, \, c^2=cb=c \big\rangle
\]
of order five; the multiplication table of this semigroup is given by
\[
\begin{array}{c|ccccc}
E  & \ 0 \, & a & ac & b  & c  \\ \hline
0  & \ 0 \, & 0 & 0  & 0  & 0  \\
a  & \ 0 \, & 0 & 0  & 0  & ac \\
ac & \ 0 \, & 0 & 0  & ac & ac \\
b  & \ 0 \, & a & ac & b  & b  \\
c  & \ 0 \, & a & ac & c  & c
\end{array}
\]
The monoid~$E^1$ was first investigated in Lee and Li~\cite[\S14]{LL11}, where its identities were shown to be finitely axiomatized by
\begin{subequations} \label{basis E}
\begin{align}
x^3 \approx x^2, \quad x^2 yx & \approx xyx, \quad xyx^2 \approx xyx, \label{basis Ea} \\
xy^2x & \approx x^2y^2. \label{basis Eb}
\end{align}
\end{subequations}
The main aim of the present section is to show that the variety
\[
\bbE^1 = \Vm\{E^1\}
\]
is of~\one.
For this purpose, the monoid~$Q^1$ obtained from the semigroup
\[
Q = \big\langle a,b,c \,\big|\, a^2 = a, \, ab = b, \, ca = c, \, ac = ba = cb = 0 \big\rangle
\]
of order five plays an important role.
It is routinely checked that the monoids~$L_2^1$, $B_0^1$, and~$Q^1$ satisfy the identities~\eqref{basis E} so the varieties $\bbL_2^1$, $\bbB_0^1$, and
\[
\bbQ^1 = \Vm\{Q^1\}
\]
are subvarieties of~$\bbE^1$.

In~\S\ref{subsec: decomposition}, the lattice $\fL(\bbE^1)$ is decomposed into the union of the lattice ${ \fL(\bbL_2^1 \vee \bbQ^1) }$ and the interval ${ [\bbL_2^1 \vee \bbQ^1,\bbE^1] }$.
The lattice ${ \fL(\bbL_2^1 \vee \bbQ^1) }$ is shown in~\S\ref{subsec: lower lattice} to be finite, while the interval ${ [\bbL_2^1 \vee \bbQ^1,\bbE^1] }$ is shown in~\S\ref{subsec: interval} to be countably infinite.
Therefore the monoid~$E^1$ generates a monoid variety of~\one.
A complete description of the lattice $\fL(\bbE^1)$ is given in~\S\ref{subsec: lattice}.

\subsection{A decomposition of $\fL(\bbE^1)$} \label{subsec: decomposition}

\begin{proposition} \label{P: decomposition}
${ \fL(\bbE^1) = \fL(\bbL_2^1 \vee \bbQ^1) \cup [\bbL_2^1 \vee \bbQ^1,\bbE^1] }$\up.
\end{proposition}

\begin{lemma}[Lee~{\cite[\S4]{Lee12}}] \label{L: LB}
The variety ${ \bbL_2^1 \vee \bbB_0^1 }$ is defined by the identities~\eqref{basis E} and $xyxzx \approx xyzx$\up.
\end{lemma}

\begin{lemma}[Lee~{\cite[\S5]{Lee14b}}] \label{L: Q} \quad
\begin{enumerate}
\item[\rm(i)] The variety~$\bbQ^1$ is Cross\up.
\item[\rm(ii)] The variety~$\bbQ^1$ is defined by the identities~\eqref{basis Ea} and $x^2y^2 \approx y^2x^2$\up.
\item[\rm(iii)] The lattice $\fL(\bbQ^1)$ is given in Figure~\ref{F: lattices}\up.
\end{enumerate}
\end{lemma}

\begin{figure}[htbp]
\begin{picture}(70,130)(0,0) \setlength{\unitlength}{0.55mm}
\put(20,70){\circle*{2}}
\put(20,60){\circle*{2}}
\put(10,50){\circle*{2}} \put(30,50){\circle*{2}}
\put(20,40){\circle*{2}}
\put(20,30){\circle*{2}}
\put(20,20){\circle*{2}}
\put(20,10){\circle*{2}}
\put(20,70){\line(0,-1){10}}
\put(20,60){\line(-1,-1){10}} \put(20,60){\line(1,-1){10}}
\put(10,50){\line(1,-1){10}} \put(30,50){\line(-1,-1){10}}
\put(20,40){\line(0,-1){30}}
\put(21,73){\makebox(0,0)[b]{$\bbQ^1$}}
\put(24,62){\makebox(0,0)[l]{$\bbB_0^1$}}
\put(08,50){\makebox(0,0)[r]{$\bbI^1$}} \put(33,50){\makebox(0,0)[l]{$\bbJ^1$}}
\put(23,41){\makebox(0,0)[tl]{$\bbM\{xy\}$}}
\put(23,29){\makebox(0,0)[l]{$\bbM\{x\}$}}
\put(23,19){\makebox(0,0)[l]{$\bbM\varnothing$}}
\put(20,06){\makebox(0,0)[t]{$\mathbf{0}$}}
\end{picture}\caption{The lattice $\fL(\bbQ^1)$}
\label{F: lattices}
\end{figure}

\begin{lemma} \label{L: decomposition}
Let~$\bbV$ be any subvariety of~$\bbE^1$\up.
Then one of the following holds\up:
\begin{enumerate}[\rm(i)]
\item ${ \bbL_2^1 \vee \bbQ^1 \subseteq \bbV }$\up;

\item ${ \bbV \subseteq \bbL_2^1 \vee \bbB_0^1 }$\up;

\item ${ \bbV \subseteq \bbQ^1 }$\up.
\end{enumerate}
\end{lemma}

\begin{proof}
Suppose that ${ \bbL_2^1 \vee \bbQ^1 \nsubseteq \bbV }$ so that either ${ L_2^1 \notin \bbV }$ or ${ Q^1 \notin \bbV }$.
Then by Lemma~\ref{L: S S1}, either ${ L_2 \notin \Vs\bbV }$ or ${ Q \notin \Vs\bbV }$.

\medskip

\paragraph{{\it Case}~1} ${ L_2 \notin \Vs\bbV }$.
Then it follows from Almeida~\cite[Proposition~10.10.2(b)]{Alm94} that the variety $\Vs\bbV$ satisfies the identity
\begin{equation}
x^2(y^2x^2)^2 \approx (y^2x^2)^2. \label{excl L2}
\end{equation}
Since
\begin{align*}
y^2x^2 & \stackrel{\eqref{basis Ea}}{\approx} (y^4x^2)x^2 \stackrel{\eqref{basis Eb}}{\approx} y^2x^2y^2x^2 \stackrel{\eqref{excl L2}}{\approx} x^2(y^2x^2y^2)x^2 \\
& \stackrel{\eqref{basis Eb}}{\approx} x^2y^4x^4 \stackrel{\eqref{basis Eb}}{\approx} x^6y^4 \stackrel{\eqref{basis Ea}}{\approx} x^2y^2,
\end{align*}
the variety $\Vs\bbV$ satisfies the identity ${ x^2y^2 \approx y^2x^2 }$.
Hence the inclusion ${ \bbV \subseteq \bbQ^1 }$ holds by Lemma~\ref{L: Q}(ii).

\medskip

\paragraph{{\it Case}~2} ${ Q \notin \Vs\bbV }$.
Then it follows from Almeida~\cite[Lemma~6.5.14]{Alm94} that the variety $\Vs\bbV$ satisfies one of the identities:
\begin{align}
x^2yx^2zx^2 & \approx x^2yzx^2, \label{excl Qa} \\
(x^2yx^2)^2 & \approx x^2yx^2. \label{excl Qb}
\end{align}
If ${ \Vs\bbV \models \eqref{excl Qa} }$, then
\[
xyxzx \stackrel{\eqref{basis Ea}}{\approx} x^2yx^2zx^2 \stackrel{\eqref{excl Qa}}{\approx} x^2yzx^2 \stackrel{\eqref{basis Ea}}{\approx} xyzx;
\]
if ${ \Vs\bbV \models \eqref{excl Qb} }$, then
\begin{align*}
xyxzx & \stackrel{\eqref{basis Ea}}{\approx} x^2yx^2zx^2 \stackrel{\eqref{excl Qb}}{\approx} (x^2yx^2zx^2)^2 \stackrel{\eqref{basis Ea}}{\approx} (xy^2xz^2x^2)^2 \\
& \stackrel{\eqref{basis Eb}}{\approx} (x^2y^2z^2x^2)^2 \stackrel{\eqref{basis Ea}}{\approx} (x^2yzx^2)^2 \stackrel{\eqref{excl Qb}}{\approx} x^2yzx^2 \stackrel{\eqref{basis Ea}}{\approx} xyzx.
\end{align*}
In any case, the variety $\Vs\bbV$ satisfies the identity ${ xyxzx \approx xyzx }$.
Therefore the inclusion ${ \bbV \subseteq \bbL_2^1 \vee \bbB_0^1 }$ holds by Lemma~\ref{L: LB}.
\end{proof}

\begin{proof}[Proof of Proposition~\ref{P: decomposition}]
The inclusion ${ \bbB_0^1 \subseteq \bbQ^1 }$ from Figure~\ref{F: lattices} implies the inclusion ${ \bbL_2^1 \vee \bbB_0^1 \subseteq \bbL_2^1 \vee \bbQ^1 }$, while the inclusion ${ \bbQ^1 \subseteq \bbL_2^1 \vee \bbQ^1 }$ is obvious.
Therefore ${ \bbL_2^1 \vee \bbB_0^1 }$ and~$\bbQ^1$ are subvarieties of ${ \bbL_2^1 \vee \bbQ^1 }$.
The result then follows from Lemma~\ref{L: decomposition}.
\end{proof}

\subsection{The lattice $\fL(\bbL_2^1 \vee \bbQ^1)$} \label{subsec: lower lattice}

\begin{proposition} \label{P: lower lattice}
The lattice ${ \fL(\bbL_2^1 \vee \bbQ^1) }$ is given in Figure~\ref{F: lower lattice}\up.
\end{proposition}

\begin{proof}
Recall from the proof of Proposition~\ref{P: decomposition} that ${ \bbL_2^1 \vee \bbB_0^1 }$ and~$\bbQ^1$ are subvarieties of ${ \bbL_2^1 \vee \bbQ^1 }$.
By Lemma~\ref{L: decomposition}, any proper subvariety of ${ \bbL_2^1 \vee \bbQ^1 }$ is contained in either ${ \bbL_2^1 \vee \bbB_0^1 }$ or~$\bbQ^1$.
Therefore ${ \bbL_2^1 \vee \bbB_0^1 }$ and~$\bbQ^1$ are the only maximal proper subvarieties of ${ \bbL_2^1 \vee \bbQ^1 }$, whence the result holds by Figures~\ref{F: LNR} and \ref{F: lattices}.
\end{proof}

\begin{figure}[htbp]
\begin{picture}(60,150)(0,0) \setlength{\unitlength}{0.55mm}
\put(30,80){\circle*{2}}
\put(20,70){\circle*{2}} \put(40,70){\circle*{2}}
\put(10,60){\circle*{2}} \put(30,60){\circle*{2}}
\put(20,50){\circle*{2}} \put(40,50){\circle*{2}}
\put(10,40){\circle*{2}} \put(30,40){\circle*{2}}
\put(30,30){\circle*{2}}
\put(30,20){\circle*{2}}
\put(30,10){\circle*{2}}
\put(30,80){\line(-1,-1){20}} \put(30,80){\line(1,-1){10}}
\put(20,70){\line(1,-1){20}} \put(40,70){\line(-1,-1){20}}
\put(10,60){\line(0,-1){20}} \put(10,60){\line(1,-1){20}}
\put(40,50){\line(-1,-1){10}}
\put(10,40){\line(1,-1){20}} \put(30,40){\line(0,-1){30}}
\put(30,83){\makebox(0,0)[b]{$\bbL_2^1 \vee \bbQ^1$}}
\put(19,72){\makebox(0,0)[br]{$\bbL_2^1 \vee \bbB_0^1$}} \put(43,70){\makebox(0,0)[l]{$\bbQ^1$}}
\put(07,62){\makebox(0,0)[r]{$\bbL_2^1 \vee \bbM(x)$}} \put(35,60){\makebox(0,0)[l]{$\bbB_0^1$}}
\put(22,47){\makebox(0,0)[tr]{$\bbI^1$}} \put(43,50){\makebox(0,0)[l]{$\bbJ^1$}}
\put(07,40){\makebox(0,0)[r]{$\bbL_2^1$}} \put(33,41){\makebox(0,0)[tl]{$\bbM(xy)$}}
\put(33,29){\makebox(0,0)[l]{$\bbM(x)$}}
\put(33,19){\makebox(0,0)[l]{$\bbM(1)$}}
\put(30,06){\makebox(0,0)[t]{$\mathbf{0}$}}
\end{picture}\caption{The lattice $\fL(\bbL_2^1 \vee \bbQ^1)$}
\label{F: lower lattice}
\end{figure}

\subsection{The interval $[{\bbL_2^1 \vee \bbQ^1},\bbE^1]$} \label{subsec: interval}

For each ${ n \geq 1 }$, define the identity
\[
\sigma_n : \bigg(\prod_{i=1}^n(\be_ih_i)\bigg)x^2y^2 \approx \bigg(\prod_{i=1}^n(\be_ih_i)\bigg)y^2x^2,
\]
where
\[
\be_i = \begin{cases}
\,x^2 & \text{if~$i$ is odd}, \\
\,y^2 & \text{if~$i$ is even}.
\end{cases}
\]
For instance, the first three identities are
\begin{align*}
\sigma_1 &: x^2h_1 \cdot x^2y^2 \approx x^2h_1 \cdot y^2x^2, \\
\sigma_2 &: x^2h_1 \cdot y^2h_2 \cdot x^2y^2 \approx x^2h_1 \cdot y^2h_2 \cdot y^2x^2, \\
\sigma_3 &: x^2h_1 \cdot y^2h_2 \cdot x^2h_3 \cdot x^2y^2 \approx x^2h_1 \cdot y^2h_2 \cdot x^2h_3 \cdot y^2x^2.
\end{align*}
Let $\sigma_\infty$ denote the identity ${ x^2y^2hx^2y^2 \approx x^2y^2hy^2x^2 }$.

\begin{proposition} \label{P: chain}
The varieties in the interval ${ [\bbL_2^1 \vee \bbQ^1,\bbE^1] }$ constitute the chain
\begin{equation}
\bbE^1\{\sigma_1\} \subset \bbE^1\{\sigma_2\} \subset \cdots \subset \bbE^1\{\sigma_\infty\} \subset \bbE^1. \label{D: inclusions}
\end{equation}
Consequently\up, ${ \bbL_2^1 \vee \bbQ^1 = \bbE^1\{\sigma_1\} }$\up.
\end{proposition}

The proof of Proposition~\ref{P: chain} is given at the end of the subsection.

\begin{lemma} \label{L: order}
The inclusions in~\eqref{D: inclusions} hold\up.
\end{lemma}

\begin{proof}
There are three cases.

\medskip

\paragraph{{\it Case}~1} ${ \bbE^1\{\sigma_n\} \subset \bbE^1\{\sigma_{n+1} \} }$ for any finite~$n$.
The inclusion ${ \bbE^1\{\sigma_n\} \subseteq \bbE^1\{\sigma_{n+1} \} }$ holds because~$\sigma_{n+1}$ is obtained by performing the substitution ${ h_n \mapsto h_n\be_{n+1}h_{n+1} }$ in~$\sigma_n$.
Further, it is easily shown that by applying the identities ${ \{ \eqref{basis E},\sigma_{n+1} \} }$ to the word ${ \big(\prod_{i=1}^n(\be_ih_i)\big)x^2y^2 }$ on the left side of~$\sigma_n$, a word of the form
\[
\bigg(\prod_{i=1}^n(\be_i'h_i)\bigg)\ba\bb\bc
\]
is obtained, where ${ \be_1',\be_3',\ldots,\ba \in \{x\}^+ }$, ${ \be_2',\be_4',\ldots,\bb \in \{y\}^+ }$, and ${ \bc \in \{x,y\}^* }$.
It follows that~$\sigma_n$ is not deducible from ${ \{ \eqref{basis E},\sigma_{n+1} \} }$, whence ${ \bbE^1\{\sigma_n\} \neq \bbE^1\{\sigma_{n+1} \} }$.

\medskip

\paragraph{{\it Case}~2} ${ \bbE^1\{\sigma_n\} \subset \bbE^1\{\sigma_\infty \} }$ for any finite~$n$.
In view of part~(i), it suffices to assume that ${ n=2r }$ is even.
Since
\[
x^2y^2 h x^2y^2 \stackrel{\eqref{basis E}}{\approx} (x^2xy^2y)^r h x^2y^2 \stackrel{\sigma_n}{\approx} (x^2xy^2y)^r hy^2x^2 \stackrel{\eqref{basis E}}{\approx} x^2y^2 h y^2x^2,
\]
the inclusion ${ \bbE^1\{\sigma_n\} \subseteq \bbE^1\{\sigma_\infty \} }$ holds.
The nonequality ${ \bbE^1\{\sigma_n\} \neq \bbE^1\{\sigma_\infty \} }$ follows from an argument similar to part~(i).

\medskip

\paragraph{{\it Case}~3} ${ \bbE^1\{\sigma_\infty\} \subset \bbE^1 }$.
The inclusion ${ \bbE^1\{\sigma_\infty\} \subseteq \bbE^1 }$ vacuously holds.
But since
\[
b^2c^2ab^2c^2 = 0 \neq ac = b^2c^2ac^2b^2
\]
in~$E^1$, the identity~$\sigma_\infty$ does not hold in~$E^1$. Hence ${ \bbE^1\{\sigma_\infty\} \neq \bbE^1 }$.
\end{proof}

It remains to verify that the varieties in the interval ${ [ \bbL_2^1 \vee \bbQ^1,\bbE^1] }$ are precisely those in~\eqref{D: inclusions}.
The remainder of this subsection is devoted to this task.

A word of the form ${ x_1^2 x_2^2 \cdots x_n^2 }$, where ${ n \geq 0 }$ and ${ x_1,x_2,\ldots,x_n }$ are distinct variables, is called a \textit{product of distinct squares}.
Note that by definition, the empty word is also a product of distinct squares.
A word~$\bu$ is in \textit{canonical form} if
\[
\bu = \bu_0 \prod_{i=1}^n (h_i\bu_i)
\]
for some ${ n \geq 0 }$, where the variables ${ h_1,h_2,\ldots,h_n }$ are simple in~$\bu$ and the possibly empty words ${ \bu_0,\bu_1,\ldots,\bu_n }$ are products of distinct squares.
It is convenient to call~$\bu_i$ the \textit{$i$th block} of~$\bu$.
Note that if ${ n=0 }$, then ${ \bu = \bu_0 }$ does not contain any simple variable.

\begin{lemma}[Lee and Li~{\cite[Lemma~14.2]{LL11}}] \label{L: canonical}
Given any word~$\bu$\up, there exists a word~$\bu'$ in canonical form such that the identity ${ \bu \approx \bu' }$ is deducible from~\eqref{basis E}\up.
\end{lemma}

\begin{lemma} \label{L: words L Q} \quad
\begin{enumerate}
\item[\rm(i)] Suppose that
\[
\bu = \bu_0 \prod_{i=1}^n (h_i\bu_i) \quad \text{and} \quad \bv = \bv_0 \prod_{i=1}^{n'} (h_i'\bv_i)
\]
are words in canonical form\up.
Then ${ \bbQ^1 \models \bu \approx \bv }$ if and only if ${n=n'}$\up, ${h_i=h_i'}$ for all~$i$\up, and ${\con(\bu_i)=\con(\bv_i)}$ for all~$i$\up.

\item[\rm(ii)] The variety~$\bbL_2^1$ satisfies an identity ${ \bu \approx \bv }$ if and only if ${ \ini(\bu) = \ini(\bv) }$\up.
\end{enumerate}
\end{lemma}

\begin{proof}
Part~(i) can be extracted from Lee and Li~\cite[Proof of Proposition~4.3]{LL11}, while Part~(ii) is well known and easily verified.
\end{proof}

\begin{lemma} \label{L: canonical int}
Suppose that~$\bu$ and~$\bv$ are any distinct words in canonical form such that ${ \bbL_2^1 \vee \bbQ^1 \models \bu \approx \bv }$\up.
Then ${ \bbE^1\{\bu \approx \bv\} = \bbE^1 \Lambda }$ for some finite set~$\Lambda$ of identities of the form
\begin{equation}
\bigg( \prod_{i=1}^m (\bp_ih_i)\bigg) x^2 y^2 \approx \bigg( \prod_{i=1}^m (\bp_ih_i)\bigg) y^2x^2, \label{D: Lambda}
\end{equation}
where ${ \bp_i \in \{ 1,x^2,y^2,x^2y^2,y^2x^2 \} }$ for all~$i$ with ${ (\bp_1,\bp_2,\ldots,\bp_m) \neq (1,1,\ldots,1) }$ and ${ m \geq 1 }$\up.
\end{lemma}

\begin{proof}
Since ${ \bbQ^1 \models \bu \approx \bv }$ where the words~$\bu$ and~$\bv$ are in canonical form, it follows from Lemma~\ref{L: words L Q}(i) that
\[
\bu = \bu_0 \prod_{i=1}^n (h_i\bu_i) \quad \text{and} \quad \bv = \bv_0 \prod_{i=1}^n (h_i\bv_i)
\]
for some ${ n \geq 0 }$, where ${ h_1,h_2,\ldots,h_n }$ are simple variables and
\begin{enumerate}
\item[(a)] ${ \bu_0,\bu_1, \ldots, \bu_n,\bv_0,\bv_1, \ldots, \bv_n }$ are products of distinct squares
\end{enumerate}
such that
\begin{enumerate}
\item[(b)] ${ \con(\bu_i) = \con(\bv_i) }$ for all~$i$.
\end{enumerate}
Further, the assumption ${ \bbL_2^1 \models \bu \approx \bv }$ and Lemma~\ref{L: words L Q}(ii) imply that
\begin{enumerate}
\item[(c)] ${ \ini(\bu) = \ini(\bv) }$.
\end{enumerate}
In particular, ${ \bu_0 = \bv_0 }$ by (a)--(c).
Since ${ \bu \neq \bv }$ by assumption, it follows that ${ n \geq 1 }$, and there exists a least ${ \ell \in \{1,2,\ldots,n \} }$ such that~$\bu$ and~$\bv$ do not share the same $\ell$th block, that is, ${ \bu_i=\bv_i }$ for all ${ i<\ell }$ while ${\bu_\ell\neq\bv_\ell }$.
For the remainder of this proof, it is shown that
\begin{enumerate}
\item[(\dag)] there exist some word~$\bv'$ in canonical form sharing the same $i$th block with~$\bu$ for all ${ i \leq \ell }$ and some finite set~$\Lambda$ of identities from~\eqref{D: Lambda} such that \[\bbE^1\{ \bu \approx \bv\} = \bbE^1 (\{ \bu \approx \bv'\} \cup \Lambda).\]
\end{enumerate}
The procedure used to establish~(\dag) can be repeated to complete the proof of the present lemma.

Let~$\bq$ be the longest suffix shared by~$\bu_\ell$ and~$\bv_\ell$.
Then~$\bu$ and~$\bv$ are of the form
\[
\bu = \ba \cdot {h_\ell \underbrace{\,\cdots x^2 \bq\,}_{\bu_\ell}} \cdot \bb \quad \text{and} \quad \bv = \ba \cdot {h_\ell \underbrace{\,\cdots x^2 y_1^2y_2^2 \cdots y_k^2 \bq\,}_{\bv_\ell}} \cdot \bc,
\]
where
\[
\ba = \bu_0 \prod_{i=1}^{\ell-1} (h_i\bu_i) = \bv_0 \prod_{i=1}^{\ell-1} (h_i\bv_i), \quad \bb = \prod_{i=\ell+1}^n (h_i\bu_i), \quad \text{and} \quad \bc = \prod_{i=\ell+1}^n (h_i\bv_i).
\]
Since ${ \bv_\ell = \cdots x^2 y_1^2y_2^2 \cdots y_k^2 \bq }$ is a product of distinct squares, ${ y_1 \notin \con(\bq) }$.
Further, since ${ \con(\bu_\ell) = \con(\bv_\ell) }$ by~(b), it follows that
\[
\bu[x,y_1] = \ba[x,y_1] \cdot y_1^2x^2 \cdot \bb[x,y_1] \quad \text{and} \quad \bv[x,y_1] = \ba[x,y_1] \cdot x^2y_1^2 \cdot \bc[x,y_1].
\]
Now if ${ \ba[x,y_1] = 1 }$, then ${ \bu[x,y_1] = y_1^2x^2 \cdots }$ and ${ \bv[x,y_1] = x^2y_1^2 \cdots }$; but this contradicts~(c).
Therefore ${ \ba[x,y_1] \neq 1 }$, whence
\begin{enumerate}
\item[(e)] the words ${ \bu_i[x,y_1]=\bv_i[x,y_1] }$, where ${ i \in \{ 0,1,\ldots, \ell-1 \} }$, are not all empty.
\end{enumerate}
Since
\begin{align*}
\bu[x,y_1,h_1,h_2,\ldots,h_\ell] & \makebox[0.32in]{$=$} \bv_0[x,y_1] \bigg( \prod_{i=1}^{\ell-1} \big(h_i \bv_i[x,y_1]\big)\bigg) h_\ell y_1^2 x^2 \bb[x,y_1] \\
& \makebox[0.32in]{$\stackrel{\eqref{basis E}}{\approx}$} \bv_0[x,y_1] \bigg( \prod_{i=1}^{\ell-1} \big(h_i \bv_i[x,y_1]\big)\bigg) h_\ell y_1^2 x^2
\end{align*}
and
\begin{align*}
\bv[x,y_1,h_1,h_2,\ldots,h_\ell] & \makebox[0.32in]{$=$} \bv_0[x,y_1] \bigg( \prod_{i=1}^{\ell-1} \big(h_i \bv_i[x,y_1]\big)\bigg) h_\ell x^2 y_1^2 \bc[x,y_1] \\
& \makebox[0.32in]{$\stackrel{\eqref{basis E}}{\approx}$} \bv_0[x,y_1] \bigg( \prod_{i=1}^{\ell-1} \big(h_i \bv_i[x,y_1]\big)\bigg) h_\ell x^2 y_1^2,
\end{align*}
the identities ${ \{ \eqref{basis E}, \bu \approx \bv \} }$ imply the identity
\[
\lambda_1 : \bv_0[x,y_1] \bigg( \prod_{i=1}^{\ell-1} \big(h_i \bv_i[x,y_1]\big)\bigg) h_\ell x^2 y_1^2 \approx \bv_0[x,y_1] \bigg( \prod_{i=1}^{\ell-1} \big(h_i \bv_0[x,y_1]\big)\bigg) h_\ell y_1^2 x^2,
\]
which, by~(e), is from~\eqref{D: Lambda}.
Hence ${ \bbE^1 \{\bu \approx \bv\} = \bbE^1\{ \bu \approx \bv, \lambda_1 \} }$.

Now since~$\bv_i$ is a product of distinct squares, ${ \bv_i[x,y_1] \in \{ 1,x^2,y_1^2,x^2y_1^2,y_1^2x^2 \} }$.
It is thus routinely checked that ${ \eqref{basis E} \vdash \bv_i \approx  \bv_i \cdot (\bv_i[x,y_1]) }$.
Since
\begin{align*}
\bv & \makebox[0.32in]{$=$} \bv_0 \bigg(\prod_{i=1}^{\ell-1} (h_i\bv_i)\bigg) {h_\ell \overbrace{\cdots x^2 y_1^2y_2^2 \cdots y_k^2 \bq}^{\bv_\ell}} \cdot \bc \\
& \makebox[0.32in]{$\stackrel{\eqref{basis E}}{\approx}$} \bv_0 \cdot (\bv_0[x,y_1]) \bigg(\prod_{i=1}^{\ell-1} \big(h_i\bv_i \cdot (\bv_i[x,y_1])\big)\bigg) h_\ell \cdots x^2 y_1^2y_2^2 \cdots y_k^2 \bq \cdot \bc \\
& \makebox[0.32in]{$\stackrel{\lambda_1}{\approx}$} \bv_0 \cdot (\bv_0[x,y_1]) \bigg(\prod_{i=1}^{\ell-1} \big(h_i\bv_i \cdot (\bv_i[x,y_1])\big)\bigg) h_\ell \cdots y_1^2x^2 y_2^2 \cdots y_k^2 \bq \cdot \bc \\
& \makebox[0.32in]{$\stackrel{\eqref{basis E}}{\approx}$} \bv_0 \bigg(\prod_{i=1}^{\ell-1} (h_i\bv_i)\bigg) h_\ell \cdots y_1^2x^2 y_2^2 \cdots y_k^2 \bq \cdot \bc = \bv^{(1)},
\end{align*}
it follows that
\[
\bbE^1 \{\bu \approx \bv\} = \bbE^1\{ \bu \approx \bv^{(1)}, \lambda_1\},
\]
where~$\bv^{(1)}$ is a word in canonical form that is obtained from~$\bv$ by interchanging the factors~$x^2$ and~$y_1^2$.
The same argument can be repeated on the identity ${ \bu \approx \bv^{(1)} }$, resulting in
\[
\bbE^1 \{\bu \approx \bv^{(1)}\} = \bbE^1\{ \bu \approx \bv^{(2)}, \lambda_2\},
\]
where~$\lambda_2$ is the identity
\[
\bv_0[x,y_2] \bigg( \prod_{i=1}^{\ell-1} \big(h_i \bv_i[x,y_2]\big)\bigg) h_\ell x^2 y_2^2 \approx \bv_0[x,y_2] \bigg( \prod_{i=1}^{\ell-1} \big(h_i \bv_0[x,y_2]\big)\bigg) h_\ell y_2^2 x^2
\]
from~\eqref{D: Lambda} and
\[
\bv^{(2)} = \bv_0 \bigg(\prod_{i=1}^{\ell-1} (h_i\bv_i)\bigg) h_\ell \cdots y_1^2y_2^2x^2 y_3^2 \cdots y_k^2 \bq \cdot \bc
\]
is a word in canonical form that is obtained from~$\bv^{(1)}$ by interchanging the factors~$x^2$ and~$y_2^2$.
Continuing in this manner,
\[
\bbE^1 \{\bu \approx \bv^{(k-1)}\} = \bbE^1\{ \bu \approx \bv^{(k)}, \lambda_k\},
\]
where~$\lambda_k$ is the identity
\[
\bv_0[x,y_k] \bigg( \prod_{i=1}^{\ell-1} \big(h_i \bv_i[x,y_k]\big)\bigg) h_\ell x^2 y_k^2 \approx \bv_0[x,y_k] \bigg( \prod_{i=1}^{\ell-1} \big(h_i \bv_0[x,y_k]\big)\bigg) h_\ell y_k^2 x^2
\]
from~\eqref{D: Lambda} and
\[
\bv^{(k)} = \bv_0 \bigg(\prod_{i=1}^{\ell-1} (h_i\bv_i)\bigg) h_\ell \cdots y_1^2y_2^2 \cdots y_k^2x^2 \bq \cdot \bc
\]
is a word in canonical form that is obtained from~$\bv^{(k-1)}$ by interchanging the factors~$x^2$ and~$y_k^2$.
Consequently,
\begin{align*}
\bbE^1 \{\bu \approx \bv\} & \makebox[0.2in]{$=$} \bbE^1\{ \bu \approx \bv^{(1)}, \lambda_1\} \\
& \makebox[0.2in]{$=$} \bbE^1\{ \bu \approx \bv^{(2)}, \lambda_1, \lambda_2\} \\
& \makebox[0.2in]{$\vdots$} \\
& \makebox[0.2in]{$=$} \bbE^1\{ \bu \approx \bv^{(k)}, \lambda_1, \lambda_2,\ldots,\lambda_k\},
\end{align*}
where the $\ell$th block of~$\bu$ and the $\ell$th block of~$\bv^{(k)}$ share the longer suffix~$x^2\bq$.

The preceding procedure can be repeated until~(\dag) is established.
\end{proof}

\medskip

\begin{proof}[Proof of Proposition~\ref{P: chain}]
Let~$\bbV$ be any variety in the interval ${ [\bbL_2^1 \vee \bbQ^1,\bbE^1] }$.
Then by Lemma~\ref{L: canonical}, there exists some set~$\Lambda$ of identities formed by words in canonical form such that ${ \bbV = \bbE^1 \Lambda }$.
By Lemma~\ref{L: canonical int}, the identities in~$\Lambda$ can be chosen from~\eqref{D: Lambda}.
In the following, it is shown that if
\[
\lambda: \bigg( \prod_{i=1}^m (\bp_ih_i)\bigg) x^2 y^2 \approx \bigg( \prod_{i=1}^m (\bp_ih_i)\bigg) y^2x^2
\]
is any identity from~\eqref{D: Lambda}, then ${ \bbE^1\{\lambda\} = \bbE^1\{\sigma_n \} }$ for some~$n$.
This is sufficient in view of Lemma~\ref{L: order}.

It is easily shown that if ${ \bp_j = 1 }$, then the identity~$\lambda'$ obtained from~$\lambda$ by removing the variable~$h_j$ is an identity from~\eqref{D: Lambda} such that ${ \bbE^1\{ \lambda\} = \bbE^1 \{ \lambda'\} }$.
Hence generality is not lost by assuming that ${ \bp_1,\bp_2,\ldots,\bp_m \in \{ x^2,y^2,x^2y^2,y^2x^2 \} }$.
Let~$\bu$ and~$\bv$ be the words on the left and right sides of the identity~$\lambda$.

\medskip

\paragraph{{\it Case}~1} ${ \bp_1,\bp_2,\ldots,\bp_m \in \{ x^2,y^2 \} }$.
Suppose that ${ \bp_j = \bp_{j+1} }$.
By symmetry, it suffices to assume that ${ \bp_j = \bp_{j+1} = x^2 }$.
Let~$\bu'$ and~$\bv'$ be words obtained by removing the factor~$\bp_jh_j$ from~$\bu$ and~$\bv$, respectively.
Then ${ \lambda': \bu' \approx \bv' }$ is an identity from~\eqref{D: Lambda}.
Since
\begin{align*}
\bu' & \makebox[0.32in]{$=$} \bigg( \prod_{i=1}^{j-1} (\bp_ih_i)\bigg) \bp_{j+1}h_{j+1} \bigg( \prod_{i=j+2}^m (\bp_ih_i)\bigg) x^2 y^2 \\
& \makebox[0.32in]{$\stackrel{\eqref{basis E}}{\approx}$} \bigg( \prod_{i=1}^{j-1} (\bp_ih_i)\bigg) x^2x \cdot x^2h_{j+1} \bigg( \prod_{i=j+2}^m (\bp_ih_i)\bigg) x^2 y^2 \\
& \makebox[0.32in]{$\stackrel{\lambda}{\approx}$} \bigg( \prod_{i=1}^{j-1} (\bp_ih_i)\bigg) x^2x \cdot x^2h_{j+1} \bigg( \prod_{i=j+2}^m (\bp_ih_i)\bigg) y^2 x^2 \\
& \makebox[0.32in]{$\stackrel{\eqref{basis E}}{\approx}$} \bigg( \prod_{i=1}^{j-1} (\bp_ih_i)\bigg) \bp_{j+1}h_{j+1} \bigg( \prod_{i=j+2}^m (\bp_ih_i)\bigg) y^2 x^2 = \bv',
\end{align*}
the deduction ${ \eqref{basis E} \cup \{ \lambda \} \vdash \lambda' }$ holds.
The deduction ${ \eqref{basis E} \cup \{ \lambda' \} \vdash \lambda }$ is easily established.
Therefore ${ \bbE^1\{\lambda\} = \bbE^1\{\lambda'\} }$.
Hence whenever the $i$th blocks and the $(i+1)$st blocks of the words~$\bu$ and~$\bv$ are equal, then the $i$th blocks can be removed and the resulting identity will still define the same subvariety $\bbE^1\{\lambda\}$ of~$\bbE^1$.
Consequently, it can be assumed that no two consecutive ${ \bp_1, \bp_2, \ldots, \bp_m \in \{ x^2,y^2\} }$ are equal.
The identity~$\lambda$ is thus~$\sigma_n$ for some finite~$n$.

\medskip

\paragraph{{\it Case}~2} ${ \bp_j \in \{ x^2y^2,y^2x^2\} }$ for some~$j$.
By symmetry, assume that ${ \bp_j = x^2y^2 }$.
Since
\[
x^2y^2 \cdot \bu[h_m,x,y] \stackrel{\eqref{basis E}}{\approx} x^2y^2 h_m x^2y^2 \quad \text{and} \quad x^2y^2 \cdot \bv[h_m,x,y] \stackrel{\eqref{basis E}}{\approx} x^2y^2 h_m y^2x^2,
\]
the deduction ${ \eqref{basis E} \cup \{ \lambda \} \vdash \sigma_\infty }$ holds.
The deduction ${ \eqref{basis E} \cup \{ \sigma_\infty \} \vdash \lambda }$ is easily verified.
Therefore ${ \bbE^1\{\lambda\} = \bbE^1\{\sigma_\infty\} }$.
\end{proof}

\subsection{The lattice $\fL(\bbE^1)$} \label{subsec: lattice}

By Propositions~\ref{P: decomposition}, \ref{P: lower lattice}, and~\ref{P: chain}, the lattice $\fL(\bbE^1)$ is given in Figure~\ref{F: lattice}.

\begin{figure}[htbp]
\begin{picture}(100,220)(0,0) \setlength{\unitlength}{0.55mm}
\put(30,130){\circle*{2}}
\put(30,120){\circle*{2}}
\put(30,100){\circle*{2}}
\put(30,90){\circle*{2}}
\put(30,80){\circle*{2}}
\put(20,70){\circle*{2}} \put(40,70){\circle*{2}}
\put(10,60){\circle*{2}} \put(30,60){\circle*{2}}
\put(20,50){\circle*{2}} \put(40,50){\circle*{2}}
\put(10,40){\circle*{2}} \put(30,40){\circle*{2}}
\put(30,30){\circle*{2}}
\put(30,20){\circle*{2}}
\put(30,10){\circle*{2}}
\put(30,130){\line(0,-1){14}}
\put(30,80){\line(0,1){24}}
\put(30,80){\line(-1,-1){20}} \put(30,80){\line(1,-1){10}}
\put(20,70){\line(1,-1){20}} \put(40,70){\line(-1,-1){20}}
\put(10,60){\line(0,-1){20}} \put(10,60){\line(1,-1){20}}
\put(40,50){\line(-1,-1){10}}
\put(10,40){\line(1,-1){20}} \put(30,40){\line(0,-1){30}}
\put(31,134){\makebox(0,0)[b]{$\bbE^1$}}
\put(34,120){\makebox(0,0)[l]{$\bbE^1\{\sigma_\infty\}$}}
\put(30,112){\makebox(0,0){$\vdots$}}
\put(34,101){\makebox(0,0)[l]{$\bbE^1\{\sigma_3\}$}}
\put(34,91){\makebox(0,0)[l]{$\bbE^1\{\sigma_2\}$}}
\put(34.5,81){\makebox(0,0)[l]{$\bbE^1\{\sigma_1\} = \bbL_2^1 \vee \bbQ^1$}}
\put(19,72){\makebox(0,0)[br]{$\bbL_2^1 \vee \bbB_0^1$}} \put(43,70){\makebox(0,0)[l]{$\bbQ^1$}}
\put(07,62){\makebox(0,0)[r]{$\bbL_2^1 \vee \bbM\{x\}$}} \put(35,60){\makebox(0,0)[l]{$\bbB_0^1$}}
\put(22,47){\makebox(0,0)[tr]{$\bbI^1$}} \put(43,50){\makebox(0,0)[l]{$\bbJ^1$}}
\put(07,40){\makebox(0,0)[r]{$\bbL_2^1$}} \put(33,41){\makebox(0,0)[tl]{$\bbM\{xy\}$}}
\put(33,29){\makebox(0,0)[l]{$\bbM\{x\}$}}
\put(33,19){\makebox(0,0)[l]{$\bbM\varnothing$}}
\put(30,06){\makebox(0,0)[t]{$\mathbf{0}$}}
\end{picture}\caption{The lattice $\fL(\bbE^1)$}
\label{F: lattice}
\end{figure}

\begin{proposition} \quad
\begin{enumerate}
\item[\rm(i)] The variety~$\bbE^1$ is \hfb\up.

\item[\rm(ii)] Every proper subvariety of $\bbE^1\{\sigma_\infty\}$ is Cross\up.

\item[\rm(iii)] The variety $\bbE^1\{\sigma_\infty\}$ is the only non\-finitely generated subvariety of~$\bbE^1$\up.
\end{enumerate}
\end{proposition}

\begin{proof}
(i) This follows from Lemma~\ref{L: Q}(i) and Propositions~\ref{P: LNR} and~\ref{P: chain}.

(ii) Let~$\bbV$ be any proper subvariety of $\bbE^1\{\sigma_\infty\}$.
Then the variety~$\bbV$ is finitely based and small by part~(i) and Figure~\ref{F: lattice}, respectively.
Further, $\bbV$ is locally finite because the variety~$\bbE^1$ is finitely generated.
Since all locally finite, small varieties are finitely generated by Lemma~\ref{L: FG}, the variety~$\bbV$ is also finitely generated.

(iii) Since ${ \bbE^1\{\sigma_\infty\} = \bigvee_{i\geq1} \bbE^1\{\sigma_i\} }$, the variety $\bbE^1\{\sigma_\infty\}$ is non\-finitely generated by Lemma~\ref{L: FG}.
The result then follows from part~(ii).
\end{proof}

\section{Other extreme properties and open questions} \label{sec: other}

\subsection{Minimal examples}

Recall that the monoids~$A_2^1$ and~$B_2^1$ generate monoid varieties of~{\two} and the monoid~$E^1$ generates a monoid variety of~\one.
Since these monoids are of order six, it is natural to question if there exists a monoid of order less than six that generates a monoid variety of~{\two} or~\one.
By Lemma~\ref{L: min non small}, the monoid $P_2^1$ of order five is the only example.

\begin{question} \label{Q: M5}
Is the monoid variety $\Vm\{P_2^1\}$ of~{\two} or~\one?
\end{question}

Regardless of the answer to this question, $P_2^1$ is the unique minimal monoid that generates a monoid variety of either~{\two} or~\one; in the former case, $E^1$ is a minimal monoid that generates a monoid variety of~\one, while in the latter case, $A_2^1$ and~$B_2^1$ are minimal monoids that generate monoid varieties of~\two.

Kad'ourek~\cite{Kad02} has exhibited an inverse semigroup of order~165 that generates an inverse semigroup variety of~{\two}.
The smallest possible order of such an inverse semigroup is also of interest.

\subsection{Finitely universal varieties}

Following Shevrin~{et al.}~\cite{SVV09}, a variety~$\cV$ is \textit{finitely universal} if every finite lattice is order-embeddable into $\fL(\cV)$.
Vernikov and Volkov~\cite{VV98} proved that the semigroup variety~$\bH$ defined by the identity
\begin{equation}
x^2 \approx yxy \label{basis VV}
\end{equation}
is finitely universal.
It is easily shown that any identity formed by a pair of non\-simple words is deducible from~\eqref{basis VV}.
Therefore a semigroup variety contains~$\bH$, and so is finitely universal, if it possesses a basis of identities formed by non\-simple words.
Due to this, many finite semigroups generate finitely universal varieties.
For instance, the semigroups
\begin{equation}
A_0, \ A_2, \ B_0, \ B_2, \ E, \ I^1, \ J^1, \ P_2^1, \, \text{ and } \, Q \label{D: FinUni}
\end{equation}
that appeared in this article generate finitely universal semigroup varieties; the semigroups $A_0$, $B_0$, $I^1$, and $J^1$ of order four, in particular, are examples of minimal order~\cite{Lee07}.
A finitely universal variety is necessarily non\-small, but it need not be of~{\two} either; for instance, the varieties generated by most of the semigroups in~\eqref{D: FinUni} are known to be {\hfb}~\cite{Lee13} and so of~\one.

Similar to the extreme properties considered earlier in the present article, the situation changes drastically when monoid varieties are considered instead of semigroup varieties.
Presently, no explicit example of finitely universal monoid variety, finitely generated or otherwise, is known.
Since any monoid that satisfies the identity~\eqref{basis VV} is trivial, the variety~$\bH$ does not come close to providing an example.
One might conjecture that the monoid variety
\[
\mathbb{H}^1 = \Vm \{ S^1 \,|\, S \in \bH \}
\]
generated by monoids obtained from semigroups in~$\bH$ is finitely universal, but the following result demonstrates that it is also futile.

\begin{proposition}
The monoid variety $\mathbb{H}^1$ contains four subvarieties\up.
\end{proposition}

\begin{proof}
It is shown that ${ \mathbb{H}^1 = \bbM\{xy\} }$.
The result then follows from Figure~\ref{F: LNR} since the monoid variety $\bbM\{xy\}$ contains four subvarieties.

It is easily shown that the subsemigroup ${ S = \RQ\{xy\} \backslash \{1\} }$ of~$\RQ\{xy\}$ satisfies the identity~\eqref{basis VV} so that ${ S \in \bH }$.
It follows that ${ \bbM\{xy\} = \Vm\{S^1\} \subseteq \mathbb{H}^1 }$.
Conversely, the identities in ${ \Gamma = \{ x^3 \approx x^2, \, xyx \approx x^2y, \, xyx \approx yx^2 \} }$ constitute a basis for the variety $\Vs\big\{\RQ\{xy\}\big\}$; see Jackson~\cite[Lemma~4.5(ii)]{Jac05b}.
Since each identity in~$\Gamma$ is formed by a pair of non\-simple words, the inclusion ${ \bH \subseteq \Vs\big\{\RQ\{xy\}\big\} }$ holds.
By Lemma~\ref{L: S S1}, the inclusion ${ \mathbb{H}^1 \subseteq \bbM\{xy\} }$ also holds.
\end{proof}

\begin{question} \label{Q: FinUni}
Do finitely universal monoid varieties exist?
\end{question}

\subsection{Cardinalities and independent systems}

To show that a variety~$\cV$ has continuum many subvarieties, the following method is employed: a system~$\Sigma$ of identities extending the equational theory of~$\cV$ is identified and shown to be independent in the sense that each subset of~$\Sigma$ defines a distinct subvariety of~$\cV$.
As far as the authors are aware, all known examples of varieties with uncountably many subvarieties have independent identity systems in this manner; see, for example, Dolinka~\cite{Dol08}, Jackson~\cite{Jac00,Jac05a}, Jackson and McKenzie~\cite{JM06}, Je\v{z}ek~\cite{Jez76}, Kad'ourek~\cite{Kad00}, Oates-Williams~\cite{O84}, Skokov and Vernikov~\cite{SV10}, Trahtman~\cite{Tra88}, Vaughan-Lee~\cite{Vau70}, and Zhang and Luo~\cite{ZL08}.
The following questions seem to be of interest for general varieties, but also under imposed restrictions such as ``semigroup'', ``locally finite'' or ``finitely generated''.

\begin{question} \label{Q: Size} \quad
\begin{enumerate}
\item[(i)] Does every variety with uncountably many subvarieties have an independent system extending its equational theory?
\item[(ii)] Does the subvariety lattice of every variety have cardinality either $2^{\aleph_0}$ or at most~$\aleph_0$? (This is trivial unless $2^{\aleph_0} \neq \aleph_1$ is assumed.)
\end{enumerate}
\end{question}

Note that the dual of any algebraic lattice with countably many compact elements arises as an interval in some subvariety lattice~\cite[Theorem~3]{Jez76}.

\section*{Acknowledgment}

The authors thank the reviewer for a number of suggestions and for bringing Question~\ref{Q: Sapir} to their attention.

\bibliographystyle{amsplain}

\end{document}